\documentclass[12pt]{amsart}
\usepackage[colorlinks=true]{hyperref}
\usepackage[utf8]{inputenc}
\usepackage{amsmath,amsthm,amssymb,url,mathdots}
\usepackage[margin = 2.5cm]{geometry}
\usepackage{tikz}
\usepackage{ytableau}
\usepackage{color}
\usepackage{MnSymbol}
\usepackage{graphicx}
\usepackage{scalerel}
\usepackage{enumitem}
\usepackage[normalem]{ulem}
\newtheorem*{theorem}{Theorem}

\newtheorem{prop}{Proposition}

\newtheorem{definition}[prop]{Definition}
\newtheorem{example}[prop]{Example}
\newtheorem{remark}[prop]{Remark}
\newtheorem{problem}{Problem}

\renewcommand{\u}[1]{\underline{#1}}

\newcommand{\q}[1]{\mathbf #1}

\newcommand{\mtk}{{\MT(\q k)}}
\DeclareMathOperator{\id}{id}

\DeclareMathOperator{\GT}{\underline{GT}}
\DeclareMathOperator{\MT}{\underline{MT}}
\DeclareMathOperator{\SGT}{\underline{SGT}}
\DeclareMathOperator{\AP}{\underline{AP}}
\DeclareMathOperator{\AR}{\underline{AR}}
\DeclareMathOperator{\sign}{sign}

\newcommand{\Z}{{\mathbb Z}}
\newcommand{\N}{{\mathbb N}}
\newcommand{\e}{{\operatorname{E}}}
\newcommand{\seswarrow}{\swarrow \!\!\!\!\!\;\!\! \searrow}
\newcommand{\nenwarrow}{\nwarrow \!\!\!\!\!\;\!\! \nearrow}
\usepackage{enumitem}
\newcommand{\si}[2]{\u{[#1,#2]}}
\newcommand{\wt}{\widetilde}

\newcommand{\comment}[1]{}

\newcommand{\mtthree}[6]{\begin{array}{c}
#1 \\ #2 \: #3 \\ #4 \: #5 \: #6
\end{array}}

\newcommand{\mttwo}[3]{\begin{array}{c}
#1 \\ #2 \: #3
\end{array}}

\title[A bijective proof of the ASM theorem, Part I: the operator formula]{A bijective proof of the ASM theorem\\ Part I: the operator formula}

\author{Ilse Fischer}
\address{Faculty of Mathematics, University of Vienna, Austria}
\urladdr{https://www.mat.univie.ac.at/~ifischer/}
\author{Matja\v z Konvalinka}%
\address{Faculty of Mathematics and Physics, University of Ljubljana, and Institute of Mathematics, Physics and Mechanics, Slovenia}
\urladdr{http://www.fmf.uni-lj.si/~konvalinka/}
\thanks{The first author acknowledges the financial support from the Austrian Science Foundation FWF (SFB F50) and the second author acknowledges the financial support from the Slovenian Research Agency (research core funding No. P1-0294).}

\date{\today}

\begin{document}

\begin{abstract}
Alternating sign matrices are known to be equinumerous with descending plane partitions, totally symmetric self-complementary plane partitions and alternating sign tri\-angles, but no bijective proof for any of these equivalences has been found so far. In this
paper we provide the first bijective proof of the operator formula for monotone triangles, which has been the main tool for several non-combinatorial proofs of such equivalences. In this proof, signed sets and sijections (signed bijections) play a fundamental role.
\end{abstract}

\maketitle

\section{Introduction}

An \emph{alternating sign matrix} (ASM) is a square matrix with entries in $\{0,1,-1\}$ such that in each row and each column the non-zero entries alternate and sum up to $1$. Robbins and Rumsey introduced alternating sign matrices in the 1980s \cite{lambda}  when studying their \emph{$\lambda$-determinant} (a generalization of the classical determinant) and showing that the $\lambda$-deter\-mi\-nant can be expressed as a sum over all alternating sign matrices of fixed size. The classical determinant is obtained from this by setting $\lambda=-1$, in which case the sum reduces so that it extends only over all ASMs \emph{without} $-1$'s, i.e., permutation matrices, and the well-known formula of Leibniz is recovered.
Numerical experiments led Robbins and Rumsey to conjecture that the number of $n \times n$ alternating sign matrices is given by the surprisingly simple product formula
\begin{equation}
\label{asm}
\prod_{i=0}^{n-1} \frac{(3i+1)!}{(n+i)!}.
\end{equation}

\medskip

Back then the surprise was even bigger when they learned from Stanley (see \cite{BrePro99,Bre99}) that this product formula had recently also appeared in Andrews' paper \cite{And79} on his proof of the weak Macdonald conjecture, which in turn provides a formula for the number of \emph{cyclically symmetric plane partitions}. As a byproduct, Andrews had introduced \emph{descending plane partitions} and had proven that the number of descending plane partitions (DPPs) with parts at most $n$ is also equal to \eqref{asm}. Since then the problem of finding an explicit bijection between alternating sign matrices and descending plane partitions has attracted considerable attention from combinatorialists and to many of them it is a miracle that such a bijection has not been found so far. All the more so because Mills, Robbins and Rumsey had also introduced several ``statistics'' on alternating sign matrices and on descending plane partitions for which they had strong numerical evidence that the joint distributions coincide as well, see \cite{MilRobRum83}.

\medskip

There were a few further surprises yet to come. Robbins introduced a new operation on plane partitions, \emph{complementation}, and had strong numerical evidence that totally symmetric self-complementary plane partitions (TSSCPPs) in a $2n \times 2n \times 2n$-box are also counted by \eqref{asm}. Again this was further supported by statistics that have the same joint distribution as well as certain refinements, see \cite{MilRobRum86,Kra96,krattsurvey,bianecheballah}. We still lack an explicit bijection between TSSCPPs and ASMs, as well as between TSSCPPs and DPPs.

\medskip

In his collection of bijective proof problems (which is available from his webpage) Stanley says the following about the problem of finding all these bijections: ``\emph{This is one of the most intriguing open problems in the area of bijective proofs.}'' In Krattenthaler's survey on plane partitions \cite{krattsurvey} he expresses his opinion by saying: ``\emph{The greatest, still unsolved, mystery concerns the question of what plane partitions have to do with alternating sign matrices.}''

\medskip

Many of the above mentioned conjectures have since been proved by non-bijective means: Zeilberger \cite{Zei96a} was the first who proved that $n \times n$ ASMs are counted by \eqref{asm}. Kuperberg gave another shorter proof \cite{Kup96} based on the remarkable observation that the \emph{six-vertex model} (which had been introduced by physicists several decades earlier) with domain wall boundary conditions is equivalent to ASMs, see \cite{ElkKupLarPro92a,ElkKupLarPro92b}, and he used the techniques that had been developed by physicists to study this model. Andrews enumerated TSSCPPs in \cite{And94}. The equivalence of certain statistics for ASMs and of certain statistics for DPPs has been proved in \cite{BehDifZin12,BehDifZin13}, while for ASMs and TSSCPPs see \cite{Zei96b,FonZin08}, and note in particular that already in Zeilberger's first ASM paper \cite{Zei96a} he could deal with an important refinement.
Further work including the study of \emph{symmetry classes} has been accomplished; for a more detailed description of this we defer to \cite{BehFisKon17}. Then, in very recent work, alternating sign triangles (ASTs) were introduced in \cite{AyyBehFis16}, which establishes a fourth class of objects that are equinumerous with ASMs, and also in this case nobody has so far been able to construct a bijection.

\medskip

Another aspect that should be mentioned here is Okada's work \cite{Oka06} (see also \cite{stro}), which hints at a connection between ASMs and representation theory that has not yet been well understood. He observed that a certain multivariate generating function (a specialization at a root of unity of the partition function that had been introduced by physicists in their study of the six-vertex model) can be expressed---up to a power of $3$---by a single \emph{Schur polynomial}. Since Schur polynomials are generating functions of semistandard tableaux, this establishes yet another challenging open problem for combinatorialists inclined to find bijections.

\medskip

The proofs of the results briefly reviewed above contain rather long and complicated computations, and include hardly any arguments of a combinatorial flavor. In fact it seems that all ASM-related identities for which there exists a bijective proofs are trivial, with the exception of the rotational invariance of fully packed loop configurations. This was proved by Wieland  \cite{wieland} bijectively and is also used in the celebrated proof of the Razumov-Stroganov (ex-)conjecture \cite{razstrog}.

\medskip

We come now to the purpose of the current paper. This is the first paper in a planned series that seeks to give the first bijective proofs of several results described so far. The seed of the idea to do so came from a brief discussion of the first author with Zeilberger on the problem of finding such bijections at the AMS-MAA Joint Mathematics Meetings 2019. Zeilberger mentioned that such bijections can be constructed from existing ``computational'' proofs, however, most likely these bijections are complicated. The authors of the current paper agree, in fact the first author gave her ``own'' proof of the ASM theorem in \cite{Fis06,Fis07,Fis16} and expressed some speculations in this direction in the final section of the last paper. It is also not implausible that a simple satisfactory bijective proof of the ASM theorem does not exist at all. Combinatorialists have failed to find such bijections for decades now, and we may start to ask ourselves why we are not rewarded for these efforts.

\medskip

This is how the authors of the current paper decided to work on converting the proof in \cite{Fis16} into a bijective proof. After having figured out how to actually convert computations and also having shaped certain useful fundamental concepts related to  \emph{signed sets} (see Section~\ref{sec:ss}), the translation of several steps became quite straightforward; some steps were quite challenging. Then a certain type of (exciting) dynamics evolved, where the combina\-to\-rial point of view led to simplifications and other modifications, and after this process the original ``computational'' proof is in fact rather difficult to recognize. For several obvious reasons, we find it essential to check all our constructions with computer code; to name one it can possibly be used to identify new equivalent statistics.

\medskip

After the above mentioned simplifications, it seems that \emph{signs} seem to be unavoidable. After all, if there would be a simple bijective proof that avoided signs, would it not also be plausible that such a proof could be converted into a simple ``computational'' proof that avoids signs? Such a proof has also not been found so far.

\medskip

In the remainder of the introduction we discuss the result that is proved bijectively in this paper, in particular we discuss why signed enumerations seem to be unavoidable from this point of view. We also sketch a few ideas informally before giving rigorous definitions and proofs later on.

\subsection*{The operator formula} 
We use the well-known correspondence between order $n \times n$ ASMs and \emph{monotone triangles} with bottom row $1,2,\ldots,n$. A \emph{monotone triangle} is a triangular array $(a_{i,j})_{1 \le j \le i \le n}$ of integers, where the elements are usually arranged as follows
\begin{equation}
\label{triangle}
\begin{array}{ccccccccccccccccc}
  &   &   &   &   &   &   &   & a_{1,1} &   &   &   &   &   &   &   & \\
  &   &   &   &   &   &   & a_{2,1} &   & a_{2,2} &   &   &   &   &   &   & \\
  &   &   &   &   &   & \dots &   & \dots &   & \dots &   &   &   &   &   & \\
  &   &   &   &   & a_{n-2,1} &   & \dots &   & \dots &   & a_{n-2,n-2} &   &   &   &   & \\
  &   &   &   & a_{n-1,1} &   & a_{n-1,2} &  &   \dots &   & \dots   &  & a_{n-1,n-1}  &   &   &   & \\
  &   &   & a_{n,1} &   & a_{n,2} &   & a_{n,3} &   & \dots &   & \dots &   & a_{n,n} &   &   &
\end{array},
\end{equation}
such that the integers increase weakly along $\nearrow$-diagonals and $\searrow$-diagonals, and increase strictly along rows, i.e.,
$a_{i,j} \le a_{i-1,j} \le a_{i,j+1}$ and $a_{i,j} < a_{i,j+1}$  for all $i,j$ with $1 \le j < i \le n$. In order to convert an ASM into the corresponding monotone triangle, add to each entry all the entries that are in the same column above it, and record then row by row the positions of the $1$'s, see Figure~\ref{ASM-MT} for an example.

\begin{figure}
$$
\left(
\begin{matrix}
0 & 0 & 0 & 1 & 0 & 0\\
0 & 1 & 0 & -1 & 1 & 0 \\
1 & -1 & 0 & 1 & -1 & 1 \\
0 & 1  & 0 & -1 & 1 & 0 \\
0 & 0 & 0 & 1 & 0 & 0 \\
0 & 0 & 1 & 0 & 0 & 0
\end{matrix} \right) \rightarrow
\left(
\begin{matrix}
0 & 0 & 0 & 1 & 0 & 0\\
0 & 1 & 0 & 0 & 1 & 0 \\
1 & 0 & 0 & 1 & 0 & 1 \\
1 & 1  & 0 &0 & 1 & 1 \\
1 & 1 & 0 & 1 & 1 & 1 \\
1 & 1 & 1 & 1 & 1 & 1
\end{matrix} \right) \rightarrow
\begin{array}{ccccccccccc}
  & & & & & 4  & & & & &  \\
  & & & & 2 &  & 5  & & & &  \\
  & & & 1 &  & 4 &   & 6 & & &  \\
 & & 1 &  & 2  &  & 5  &  & 6 & &  \\
 &1  &  & 2 &   & 4  &   & 5  &  & 6 &  \\
1  &  & 2 &  & 3  &   & 4   &   & 5  &  & 6
\end{array}
$$
\caption{\label{ASM-MT} ASM $\rightarrow$ partial columnsums $\rightarrow$ monotone triangle}
\end{figure}

\medskip

The following \emph{operator formula}  for the number of monotone triangles with prescribed bottom row was first proved in \cite{Fis06} (see \cite{Fis10,Fis16} for simplifications and generali\-za\-tions). Note that we allow arbitrary strictly increasing bottom rows.

\begin{theorem}
\label{operator}
Let $k_1 < k_2 < \ldots < k_n$ be a sequence of strictly increasing integers. The number of monotone triangles with bottom row $k_1,\ldots,k_n$ is
\begin{equation}
\label{operatorexpr}
\prod_{1 \le p < q \le n} \left( \e_{k_p} + \e_{k_q}^{-1} - \e_{k_p} \e_{k_q}^{-1} \right) \prod_{1 \le i < j \le n} \frac{k_j-k_i+j-i}{j-i},
\end{equation}
where $\e_x$ denotes the shift operator, i.e., $\e_x p(x) = p(x+1)$.\footnote{The formula has to be understood as follows: Take $\prod_{1 \le i < j \le n} \frac{k_j-k_i+j-i}{j-i}$ and treat the $k_i$'s as indeterminates. Apply $\prod_{1 \le p < q \le n} \left( \e_{k_p} + \e_{k_q}^{-1} - \e_{k_p} \e_{k_q}^{-1} \right)$ to this polynomial to obtain another polynomial. Only then the $k_i$'s can be specialized to the actual values.}
\end{theorem}

The purpose of this paper is to provide a bijective proof of Theorem~\ref{operator}. While the operator formula is an interesting result in its own right, it has also been the main tool for proofs of several results mentioned above. This will be reviewed in the final section of this paper along with indications for future projects on converting also these proofs into bijective proofs.

\medskip

In order to be able to construct a bijective proof of Theorem~\ref{operator}, we need to interpret \eqref{operatorexpr} combinatorially. Recall that \emph{Gelfand-Tsetlin patterns} are defined as monotone triangles with the condition on the strict increase along rows being dropped, see  \cite[p.\ 313]{Sta99} or \cite[(3)]{gelfand} for the original reference\footnote{Gelfand-Tsetlin patterns with bottom row $0 \le k_1 \le k_2 \le \ldots \le k_n$ are in an easy bijective correspondence with seminstandard tableaux of shape $(k_n,k_{n-1},\ldots,k_1)$ and entries in $\{1,2,\ldots,n\}$.}. It is well known that the number of Gelfand-Tsetlin patterns with bottom row $k_1 \le k_2 \le \ldots \le k_n$ is
\begin{equation}
\label{gelfandtsetlin}
\prod_{1 \le i < j \le n} \frac{k_j-k_i+j-i}{j-i},
\end{equation}
which is the operand in the operator formula \eqref{operatorexpr}. Expanding $\prod_{1 \le p < q \le n} \left( \e_{k_p} + \e_{k_q}^{-1} - \e_{k_p} \e_{k_q}^{-1} \right)$ into $3^{\binom{n}{2}}$ monomials in $\e^{\pm 1}_{k_1}, \e^{\pm 1}_{k_2}, \ldots, \e^{\pm 1}_{k_n}$ (keeping a copy for each multiplicity), \eqref{operatorexpr} is a signed enumeration of certain Gelfand-Tsetlin patterns, where each monomial causes a deformation of the bottom row $k_1,\ldots,k_n$. It is useful to encode these deformations by \emph{arrow patterns} as defined in Section~\ref{sec:recursion}, where we choose $\swarrow$ if we pick $\e_{k_p}$ from $\e_{k_p} + \e_{k_q}^{-1} - \e_{k_p} \e_{k_q}^{-1}$, we choose $\searrow$ if we pick $\e_{k_q}^{-1}$, while we choose $\seswarrow$ if we pick $-\e_{k_p} \e_{k_q}^{-1}$. Arranging the $\binom{n}{2}$ arrows in a triangular manner so that the arrows coming from $\e_{k_p} + \e_{k_q}^{-1} - \e_{k_p} \e_{k_q}^{-1}$ are situated in the $p$-th $\swarrow$-diagonal and the $q$-th $\searrow$-diagonal, and placing $k_1,\ldots,k_n$ in the bottom row will allow us to describe the deformation coming from a particular monomial in a convenient way. The combinatorial objects associated with \eqref{operatorexpr} then consist of a pair of such an arrow pattern and a Gelfand-Tsetlin pattern where the bottom row is a deformation of $k_1,\ldots,k_n$ as prescribed by the arrow pattern. This will lead directly to the definition of \emph{shifted Gelfand-Tsetlin patterns}.

\medskip

A sign comes from picking $- \e_{k_p} \e_{k_q}^{-1}$, but there is also a more subtle appearance. The deformation induced by the arrow pattern can cause a deformation of the increasing bottom row $k_1,k_2,\ldots,k_n$  into a sequence that is not increasing. Therefore we are in need of an extension of the combinatorial interpretation of \eqref{gelfandtsetlin} to any sequence $k_1,\ldots,k_n$ of integers. Such an interpretation was given in \cite{Fis05} and is repeated below in Section~\ref{sec:gt}.

\subsection*{Outline of the bijective proof} Given a sequence $k_1 < \ldots < k_n$, it suffices to find an injective map from the set of monotone triangles with bottom row $k_1,\ldots,k_n$ to our shifted Gelfand-Tsetlin patterns associated with $k_1,\ldots,k_n$ so that the images under this map have positive signs, along with a sign-reversing involution on the set of shifted Gelfand-Tsetlin patterns that are not the image of a monotone triangle.

\medskip

We will accomplish something more general, as we will also consider an extension of monotone triangles to all integer sequences $k_1,\ldots,k_n$, see Section~\ref{sec:recursion}, along with a sign function on these objects, and prove that the operator formula also holds in this instance. To do that, we will construct a sign-reversing involution on a subset of monotone triangles, another sign-reversing involution on a subset of shifted Gelfand-Tsetlin patterns, and a sign-\emph{preser\-ving} bijection between the remaining monotone triangles and the remaining shifted Gelfand-Tsetlin patterns. Note that this is actually equivalent to the construction of a bijection between the (disjoint) union of the ``positive'' monotone triangles and the ``negative'' shifted Gelfand-Tsetlin patterns, and the (disjoint) union of the ``negative'' monotone triangles and the ``positive'' shifted Gelfand-Tsetlin patterns. We call such maps \emph{sijections} for general signed sets.

\medskip

The actual construction here will make use of the recursion underlying monotone triangles. For a monotone triangle with bottom row $k_1,\ldots,k_n$, the eligible penultimate rows $l_1,\ldots,l_{n-1}$ are those with
$$
k_1 \le l_1 \le k_2 \le l_2 \le \ldots \le l_{n-1} \le k_n,
$$
and $l_1 < l_2 < \ldots < l_{n-1}$. This establishes a recursion that can be used to construct all monotone triangles. Phrased differently, ``at'' each $k_i$ we need to sum over all $l_{i-1},l_i$ such that $l_{i-1} \le k_i \le l_i$ and $l_{i-1} < l_{i}$.\footnote{The degenerate cases $k_1$ and $k_n$ are slightly different.} However, we can split this into the following three cases:
\begin{enumerate}
\item Consider all $l_{i-1},l_i$ with $l_{i-1} < k_i \le  l_i$.
\item Consider all $l_{i-1},l_i$ with $l_{i-1} \le  k_i < l_i$.
\item Combining (1) and (2), we have done some double counting, thus we need to subtract the intersection, i.e., all $l_{i-1},l_i$ with $l_{i-1} < k_i < l_i$.
\end{enumerate}
This can be written as a recursion. The \emph{arrow rows} in Section~\ref{sec:recursion} are used to describe this recursion: we choose $\nwarrow$ ``at'' $k_i$ if we are in Case (1), $\nearrow$ in Case (2), and $\nenwarrow$ in Case (3). Our main effort will be to show ``sijectively'' that shifted Gelfand-Tsetlin patterns also fulfill this recursion.

\subsection*{Outline of the paper} The remainder of this paper is devoted to the bijective proof of Theorem~\ref{operator} (or rather, the more general version with the increasing condition on $k_1,\ldots,k_n$ dropped). In Section~\ref{sec:ss} we lay the groundwork by defining concepts like signed sets and sijections, and we extend known concepts such as disjoint union, Cartesian product and composition for ordinary sets and bijections to signed sets and sijections. The composition of sijections will use a variation of the well-known Garsia-Milne involution principle \cite{GarsiaMilne2,GarsiaMilne1}. Many of the signed sets we will be considering are signed boxes (Cartesian products of signed intervals) or at least involve them, and we define some sijections on them in Section~\ref{sec:sb}. These sijections will be the building blocks of our bijective proof later on. In
Section~\ref{sec:gt} we introduce the extended Gelfand-Tsetlin patterns and construct some related sijections. In Section~\ref{sec:recursion}, we finally define the extended monotone triangles as well as the shifted Gelfand-Tsetlin patterns (i.e., the combinatorial interpretation of \eqref{operatorexpr}), and use all the preparation to construct the sijection between monotone triangles and shifted Gelfand-Tsetlin patterns. In the final section, we discuss further projects.

\medskip

To emphasize that we are not merely interested in the fact that two signed sets have the same size, but want to use the constructed signed bijection later on, we will be using a convention that is slightly unorthodox in our field. Instead of listing out results as lemmas and theorems with their corresponding proofs, we will be using the Problem--Construction terminology. See for instance \cite{Voevodsky} and \cite{Bauer}.

\section{Signed sets and sijections} \label{sec:ss}

\subsection*{Signed sets}

A \emph{signed set} is a pair of disjoint finite sets: $\u S = (S^+,S^-)$ with $S^+ \cap S^- = \emptyset$. Equivalently, a signed set is a finite set $S$ together with a sign function  $\sign \colon S \to \{1,-1\}$. While we will mostly avoid the use of the sign function altogether (with the exception of monotone triangles defined in Section \ref{sec:recursion}), it is useful to keep this description at the back of one's mind. Note that throughout the paper, signed sets are underlined. We will write $i \in \u S$ to mean $i \in S^+ \cup S^-$.

\medskip

The \emph{size} of a signed set $\u S$ is $|\u S| = |S^+| - |S^-|$. The \emph{opposite} signed set of $\u S$ is $- \u S = (S^-,S^+)$. We have $|{- \u S}  | = -|\u S|$. The \emph{Cartesian product} of signed sets $\u S$ and $\u T$ is
$$\u S \times \u T = (S^+ \times T^+ \cup S^- \times T^-,S^+ \times T^- \cup S^- \times T^+),$$
and we can similarly (or recursively) define the Cartesian product of a finite number of signed sets. We have
$$|\u S \times \u T| = |S^+| \cdot |T^+| + |S^-| \cdot |T^-| - |S^+| \cdot |T^-| - |S^-| \cdot |T^+| = |\u S| \cdot |\u T|.$$

The intersection of signed sets $\u S$ and $\u T$ is defined as $\u S \cap \u T = (S^+ \cap T^+, S^- \cap T^-)$, while the union $\u S \cup \u T = (S^+ \cup T^+, S^- \cup T^-)$ is only defined when $S^+ \cap T^- = S^- \cap T^+ = \emptyset$. Again, we can extend these definitions to a finite family of signed sets.

\begin{example}
 One of the crucial signed sets is the \emph{signed interval}
 $$\u{[a,b]} = \begin{cases} ([a,b],\emptyset) & \mbox{if } a \leq b \\ (\emptyset,[b+1,a-1]) & \mbox{if } a > b \end{cases}$$
 for $a,b \in \Z$, {where $[a,b]$ stands for an interval in $\mathbb{Z}$ in the usual sense}. We have $\si{b+1}{a-1} = -\si a b$ and $|\u{[a,b]}| = b - a + 1$.\\
 We will also see many \emph{signed boxes}, Cartesian products of signed intervals. Note that $S^+ = \emptyset$ or $S^- = \emptyset$ for every signed box $\u S$.
\end{example}

{Signed subsets $\u T \subseteq \u S$ are defined in an obvious manner, in particular, for $s \in \u S$, we have}
$$\u{\{s\}} = \begin{cases} (\{s\},\emptyset) & \mbox{if } s \in S^+ \\ (\emptyset,\{s\}) & \mbox{if } s \in S^- \end{cases}.$$

The \emph{disjoint union} of signed sets $\u S$ and $\u T$ is the signed set
$$\u S \sqcup \u T = (\u S \times (\{0\},\emptyset)) \cup (\u T \times (\{1\},\emptyset))$$
with elements $(s,0)$ for $s \in \u S$ and $(t,1)$ for $t \in \u T$. If $\u S$ and $\u T$ are signed sets with $(S^+ \cup S^-) \cap (T^+ \cup T^-) = \emptyset$, we can identify $\u S \cup \u T$ and $\u S \sqcup \u T$.

\medskip

More generally, we can define the disjoint union of a family of signed sets $\u S_t$, where the family is indexed with a signed set $\u T$:
$$\bigsqcup_{t \in \u T} \u S_t = \bigcup_{t \in \u T} (\u S_t \times \u{\{t\}}).$$
{We get $\bigsqcup_{t \in \u{[0,1]}} \u S_t = S_0 \sqcup S_1$. For $a,b \in \mathbb{Z}$, we may also write $\bigsqcup_{i=a}^b \u S_i$ instead of $\bigsqcup_{i \in \si a b} \u S_i$.}
{As for the size, we have $$\left| \bigsqcup_{t \in \u T} \u S_t  \right|= \sum_{t \in T} | \u{S_t} | \cdot
| \u { \{ t \} } |.$$}

\medskip

The usual properties such as {associativity $(\u S \sqcup \u T) \sqcup \u U = \u S \sqcup (\u T \sqcup \u U)$ and} {distributivity}
$(\u S \sqcup \u T) \times \u U = \u S \times \u U \sqcup \u T \times \u U$ also hold. {Strictly speaking, the $=$ sign here and sometimes later on indicates that there is an obvious and natural sign-preserving bijection between the two signed sets.} We summarize a few more basic properties that will be needed in the following and that are easy to prove.
\begin{enumerate}
\item $$\bigsqcup\limits_{\q l \in \si{a_1}{b_1} \times \ldots  \times \si{a_{n}}{b_{n}}} \u S_{l_1+c_1,\ldots,l_n+c_n}  =
\bigsqcup\limits_{\q l \in \si{a_1+c_1}{b_1+c_1} \times \ldots  \times \si{a_{n}+c_{n}}{b_{n}+c_{n}}} \u S_{l_1,\ldots,l_n}$$
\item $$
\bigsqcup_{t \in \u T} \bigsqcup_{u \in \u U} \u S_{t,u} = \bigsqcup_{(u,t) \in \u U \times \u T} 
\u S_{t,u} =    \bigsqcup_{(t,u) \in \u T \times \u U} \u S_{t,u} =
\bigsqcup_{u \in \u U} \bigsqcup_{t \in \u T} \u S_{t,u}
$$
\item
$$
\bigsqcup_{t \in \bigsqcup_{u \in \u U} \u T_u} \u S_t = \bigsqcup_{u \in \u U} \bigsqcup_{t \in \u T_u} \u S_t
$$
\item $$
- \bigsqcup_{t \in \u T} \u S_t  = \bigsqcup_{t \in \u T} - \u S_t= \bigsqcup_{t \in -\u T} \u S_t.
$$
\end{enumerate}

\medskip

\subsection*{Sijections}

The role of bijections for signed sets is played by ``signed bijections'', which we call \emph{sijections}. A sijection $\varphi$ from $\u S$ to $\u T$,
$$\varphi \colon \u S \Rightarrow \u T,$$
is an involution on the set $(S^+ \cup S^-) \sqcup (T^+ \cup T^-)$ with the property $\varphi(S^+ \sqcup T^-) = S^- \sqcup T^+$, {where $\sqcup$ refers to the disjoint union for ordinary (``unsigned'') sets.}
 It follows that also $\varphi(S^- \sqcup T^+) = S^+ \sqcup T^-$. {There is an obvious sijection $\id_{\u S} \colon \u S \Rightarrow \u S$.}
\medskip

We can think of a sijection as a collection of a sign-reversing involution on a subset of $\u S$, a sign-reversing involution on a subset of $\u T$, and a {sign-preserving} matching between the remaining elements of $\u S$ with the remaining elements of $\u T$. When $S^- = T^- = \emptyset$, the signed sets can be identified with ordinary sets, and a sijection in this case is simply a bijection.

\medskip

{A sijection is a  manifestation of the fact that two signed sets have the same size. Indeed,}
if there exists a sijection $\varphi \colon \u S \Rightarrow \u T$, we have $|S^+| + |T^-| = |S^+ \sqcup T^-| = |S^- \sqcup T^+| = |S^-| + |T^+|$ and therefore $|\u S| = |S^+| - |S^-| = |T^+| - |T^-| = |\u T|$. A sijection $\varphi \colon \u S \Rightarrow \u T$ has an inverse $\varphi^{-1} \colon \u T \Rightarrow \u S$ that we obtain by identifying $(T^+ \cup T^-) \sqcup (S^+ \cup S^-)$ with $(S^+ \cup S^-) \sqcup (T^+ \cup T^-)$.

\medskip

\comment{A sijection $\varphi \colon \u S \Rightarrow \u T$ is \emph{simple} if  $\varphi(S^+) = T^+$ and  $\varphi(S^-) = T^-$.}

\medskip

For a signed set $\u S$, there is a natural sijection $\varphi$ from $\u S \sqcup (- \u S)
$
to the empty signed set $\u \emptyset = (\emptyset,\emptyset)$. Indeed, the {involution} should be defined on {$(S^+ \times \{0\} \cup S^- \times\{1\}) \cup (S^- \times \{0\} \cup S^+ \times \{1\})$ and map $S^+ \times \{0\} \cup S^- \times\{1\}$ to $S^+ \times \{1\} \cup S^- \times \{0\}$}, and {so} we can take $\varphi((s,0),0) = ((s,1),0)$, $\varphi((s,1),0) = ((s,0),0)$. {Note that in general, a sijection from a signed set $\u S$ to $\u \emptyset$ is simply a sign-reversing involution on $\u S$, in other words, a bijection between $S^+$ and $S^-$.} \comment{In other words, a sijection $\varphi \colon \u S \Rightarrow \u \emptyset$ is equivalent to a simple sijection $\varphi \colon \u S \Rightarrow - \u S$.}

\medskip

{If we have a sijection $\varphi \colon \u S \Rightarrow \u T$, there is a natural sijection $-\varphi \colon -\u S \Rightarrow -\u T$ (as a map, it is actually precisely the same).}

\medskip

{If we have sijections $\varphi_i \colon \u S_i \Rightarrow \u T_i$ for $i=0,1$, then there is a natural sijection $\varphi \colon \u S_0 \sqcup \u S_1 \Rightarrow \u T_0 \sqcup \u T_1$.}
More interesting ways to create new sijections are described below in Proposition \ref{prop:sijections}, but we will need this in our first construction {for the special case $\u S_0 = \u T_0$ and $\varphi_0 = \id_{\u S_0}$.}

\medskip

To motivate our first result, note that if $a \leq b \leq c$ or $c < b < a$, then $\si a c = \si a b \cup \si{b+1}c= \si a b \sqcup \si{b+1}c$. Of course, this does not hold in general; for $a = 1$, $b = 8$, $c = 5$, we have $\si 1 5 = (\{1,2,3,4,5\},\emptyset)$, $(\si 1 8 \sqcup \si 9 5)^+ = (\{(1,0),(2,0),(3,0),(4,0),(5,0),(6,0),(7,0),(8,0)\}$ and $(\si 1 8 \sqcup \si 9 5)^- = \{(6,1),(7,1),(8,1)\})$. The following, however, tells us that there is in general a sijection between {$\si a c$ and
$\si a b \sqcup \si{b+1}c$.}
This map will be the crucial building block for more complicated sijections.

\begin{problem} \label{prob:alpha}
 Given $a,b,c \in \Z$, construct a sijection
 $$\alpha = \alpha_{a,b,c} \colon \si a c \Rightarrow \si a b \sqcup \si{b+1}c.$$
\end{problem}
\begin{proof}[Construction]
 For $a \leq b \leq c$ and $c < b < a$, there is nothing to prove. For, say, $a \leq c < b$, we have
 $$\si a b \sqcup \si{b+1}c = (\si a c \sqcup \si{c+1} b) \sqcup \si{b+1}c = \si a c \sqcup (\si{c+1} b \sqcup (-\si {c+1}{b}))$$
 and since there is a sijection $\si{c+1} b \sqcup (-\si {c+1}{b}) \Rightarrow \u \emptyset$, we get a sijection $\si a b \sqcup \si{b+1}c \Rightarrow \si a c$. The cases $b < a \leq c$, $b \leq c < a$, and $c < a \leq b$ are analogous. \comment{Note that in all cases, $\si a b$ and  $\si{b+1}c$ are either disjoint or one set is contained in the other.}
\end{proof}

\medskip

The following proposition describes composition, Cartesian product, and disjoint union of sijections. The composition is a variant of the well-known Garsia-Milne involution principle. All the statements are easy to prove{, and the proofs are left to the reader.}

\begin{prop} \leavevmode \label{prop:sijections}
 \begin{enumerate}
  \item (Composition) Suppose that we have sijections $\varphi \colon \u S \Rightarrow \u T$ and $\psi \colon \u T \Rightarrow \u U$. For $s \in \u S$ (resp.\ $u \in \u U$), define $\psi \circ \varphi(s)$ {(resp.\ $\psi \circ \varphi(u)$)} as the last well-defined element in the sequence $s, \varphi(s), \psi(\varphi(s)), \varphi(\psi(\varphi(s))),\ldots$ (resp.\ $u, \psi(u), \varphi(\psi(u)), \psi(\varphi(\psi(u))),\ldots$). Then $\psi \circ \varphi$ is a well-defined sijection from $\u S$ to $\u U$. \comment{If $\varphi$ and $\psi$ are simple, so is $\psi \circ \varphi$.}
  \item (Cartesian product) Suppose we have sijections $\varphi_i \colon \u S_i \Rightarrow \u T_i$, $i = 1,\ldots,k$. Then $\varphi = \varphi_1 \times \cdots \times \varphi_k$, defined by
$$\varphi(s_1,\ldots,s_k) = \begin{cases} (\varphi_1(s_1),\ldots,\varphi_k(s_k)) & \mbox{if } \varphi_i(s_i) \in \u T_i \mbox{ for } i=1,\ldots,k \\ (s_1,\ldots,s_{j-1},\varphi_j(s_j),s_{j+1},\ldots,s_k) & \mbox{if } \varphi_j(s_j) \in \u S_j, \varphi_i(s_i) \in \u T_i \mbox{ for } i < j \end{cases}$$
{if $(s_1,\ldots,s_k) \in \u S_1 \times  \dots \times \u S_k$ and
$$\varphi(t_1,\ldots,t_k) = \begin{cases} (\varphi_1(t_1),\ldots,\varphi_k(t_k)) & \mbox{if } \varphi_i(t_i) \in \u S_i \mbox{ for } i=1,\ldots,k \\ (t_1,\ldots,t_{j-1},\varphi_j(t_j),t_{j+1},\ldots,t_k) & \mbox{if } \varphi_j(t_j) \in \u T_j, \varphi_i(t_i) \in \u S_i \mbox{ for } i < j \end{cases}$$
if $(t_1,\ldots,t_k) \in \u T_1  \times \dots \times \u T_k$,}
  is a well-defined sijection from $\u S_1 \times \cdots \times \u S_k$ to $\u T_1 \times \cdots \times \u T_k$. \comment{If $\varphi_i$ are all simple sijections, so is $\varphi$.}
  \item (Disjoint union) {Suppose we have signed sets $\u T, \u{\wt T}$ and a sijection $\psi \colon \u T \Rightarrow \u{\wt T}$. Further\-more, suppose that for every $t \in \u T \sqcup \u{\wt T}$, we have a signed set $\u S_t$ and a sijection $\varphi_t \colon \u S_t \Rightarrow \u S_{\psi(t)}$ satisfying $\varphi_{\psi(t)} = \varphi_t^{-1}$.  Then $\varphi = \bigsqcup_{t \in \u T \sqcup \u{\wt T}} \varphi_t$, defined by
$$\varphi(s_t,t) = \begin{cases}
(\varphi_t(s_t),t) & \mbox{if } t \in \u T \sqcup \u {\wt T}, s_t \in \u S_t, \varphi_t(s_t) \in \u S_t \\
(\varphi_t(s_t),\psi(t)) & \mbox{if } t \in \u T \sqcup \u {\wt T}, s_t \in \u S_t, \varphi_t(s_t) \in \u {S}_{\psi(t)} \\
\end{cases}
$$
is a sijection $\bigsqcup_{t \in \u T} \u S_t \Rightarrow \bigsqcup_{t \in \u {\wt T}} \u {\wt S}_t$.} \comment{If $\psi$ and $\varphi_t$ are all simple sijections, so is $\varphi$.}
\end{enumerate}
\end{prop}

{One} important special case of Proposition \ref{prop:sijections} (3) is $\u T =  \u {\wt T}$ and $\psi = \id$. We have two sets of signed sets indexed by $\u T$, $\u S_{(t,0)} =: \u{S}^{{0}}_t$ and $\u S_{(t,1)} =: \u{S}^{{1}}_t$, and sijections $\varphi_t \colon \u{S}^{{0}}_t \Rightarrow \u{S}^{{1}}_t$. By the proposition, these sijections have a disjoint union that is a sijection $\bigsqcup_{t \in \u T} \u S^{{0}}_t \Rightarrow \bigsqcup_{t \in \u T} \u S^{{1}}_t$.

\medskip

By the proposition, the relation
$$\u S \approx \u T \iff \mbox{there exists a sijection from } \u S \mbox{ to } \u T$$
is an equivalence {relation}.

\subsection*{Elementary signed sets and normal sijections}

Often, we will be interested in disjoint unions of Cartesian products of signed intervals. An element of such a signed set is a pair, consisting of a tuple of integers and an element of the indexing signed set. Intuitively, the first one is ``more important'', as the second one serves just as an index. We formalize this notion in the following definition.

\begin{definition}
 A signed set $\u A$ is \emph{elementary of dimension $n$ and depth $0$} if its elements are in $\Z^n$. A signed set $\u A$ is \emph{elementary of dimension $n$ and depth $d$}, $d \geq 1$, if it is of the form
 $$\bigsqcup_{t \in \u T} \u S_t,$$
 where $\u T$ is a signed set, and $\u S_t$ are all signed sets of dimension $n$ and depth at most $d-1$, with the depth of at least one of them equal to $d-1$. A signed set $\u A$ is \emph{elementary of dimension $n$} if it is an elementary signed set of dimension $n$ and depth $d$ for some $d \in \N$.\\
 The \emph{projection map} on an elementary set of dimension $n$ is the map
 $$\xi \colon \u A \to \Z^n$$
 defined as follows. If the depth of $\u A$ is $0$, then $\xi$ is simply the inclusion map. Once $\xi$ is defined on elementary signed sets of depth $< d$, and the depth of $\u A$ is $d$, then $\u A = \bigsqcup_{t \in \u T} \u S_t$, where $\xi$ is defined on all $\u S_t$. Then define $\xi(s,t) = \xi(s)$ for $(s,t) \in \u A$.\\
 A sijection $\psi \colon \u T \Rightarrow \u{\wt T}$ between elementary signed sets $\u T$ and $\u{\wt T}$ of the same dimension is \emph{normal} if $\xi(\psi(t)) = \xi(t)$ for all $t \in \u T \sqcup \u{\wt T}$.
\end{definition}

Simple examples of elementary signed sets are $\si a c$, $\si a b \sqcup \si {b+1} c$ and $\si a c \sqcup (\si a b \sqcup \si {b+1} c)$. They are all of dimension $1$ and depth $0$, $1$ and $2$, respectively.\footnote{{To avoid ambiguity, we should consider signed intervals in this case to be subsets of $\Z^1$ ($1$-tuples of integers), not $\Z$. Otherwise, $\si{0}{1} \sqcup \si{2}{3} = (\{(0,0),(1,0),(2,1),(3,1)\},\emptyset)$, and this can be seen either as an elementary set of dimension $1$ and depth $1$, or as an elementary signed set of dimension $2$ and depth $0$. So the interpretation depends on the ``representation'' of the set as disjoint union. Instead, we should understand $\si{0}{1} \sqcup \si{2}{3}$ to mean $ (\{((0),0),((1),0),((2),1),((3),1)\},\emptyset)$, with dimension $1$ and depth $1$. For coding, the distinction is important, but in the paper we nevertheless think of elements of signed intervals as integers.}} It is easy to see that the sijection $\alpha_{a,b,c}$ from Problem \ref{prob:alpha} is normal.

Let us illustrate this with the example $a = 1$, $b = 5$, $c = 3$. We have $\si a c = (\{1,2,3\},\emptyset)$ and $\si a b \sqcup \si {b+1} c = (\{(1,0),(2,0),(3,0),(4,0),(5,0)\},\{(4,1),(5,1)\})$. The sijection $\alpha_{1,5,3}$ is the involution on $\si 1 3 \sqcup (\si 1 5 \sqcup \si {6} 3)$ defined by
$$(1,0) \leftrightarrow ((1,0),1), \quad (2,0) \leftrightarrow ((2,0),1), \quad (3,0) \leftrightarrow ((3,0),1),$$
$$((4,0),1) \leftrightarrow ((4,1),1), \quad ((5,0),1) \leftrightarrow ((5,1),1).$$
Since $\xi(i,0) = i$ for $i = 1,2,3$, $\xi((i,0),1) = i$ for $i = 1,2,3,4,5$ and $\xi((i,1),1) = i$ for $i=4,5$, $\alpha_{1,5,3}$ is indeed normal.

\medskip

Other examples of elementary signed sets appear in the statements of Problems \ref{prob:beta} and \ref{prob:gamma} (in both cases, they are of dimension $n-1$).

\medskip

Normality is preserved under Cartesian product, disjoint union etc. For example, the sijection
\begin{multline*}
  \si{a_1}{c_1} \times \si{a_2}{c_2} \Rightarrow \\
  \si{a_1}{b_1} \times \si{a_2}{b_2} \sqcup \si{a_1}{b_1} \times \si{b_2+1}{c_2} \sqcup \si{b_1+1}{c_1} \times \si{a_2}{b_2} \sqcup \si{b_1+1}{c_1} \times \si{b_2+1}{c_2},
\end{multline*}
obtained by using $\alpha_{a_1,b_1,c_1} \times \alpha_{a_2,b_2,c_2}$ and distributivity on disjoint unions, is normal.

\medskip

The main reason normal sijections are important is that they give a very natural special case of Proposition \ref{prop:sijections} (3). Suppose that $\u T$ and $\u{\wt T}$ are elementary signed sets of dimension $n$, and that $\psi \colon \u T \Rightarrow \u{\wt T}$ is a normal sijection. Furthermore, suppose that we have a signed set $\u S_{\q k}$ for every $\q k \in \Z^n$. Then we have a sijection
$$\bigsqcup_{t \in \u T} \u S_{\xi(t)} \Rightarrow \bigsqcup_{t \in \u{\wt T}} \u S_{\xi(t)}.$$
Indeed, Proposition \ref{prop:sijections} gives us a sijection provided that we have a sijection $\varphi_t \colon \u S_{\xi(t)} \Rightarrow \u S_{\xi(\psi(t))}$ satisfying $\varphi_{\psi(t)} = \varphi_t^{-1}$ for every $t \in \u T \sqcup \u{\wt T}$. But since $\xi(\psi(t)) = \xi(t)$, we can take $\varphi_t$ to be the identity.

\comment{The following observation is crucial: Suppose $\u T \approx \u {\wt T}$, then in general we do not have
\begin{equation}
\label{indexset_equivalence}
\bigsqcup_{t \in \u T} S_t \approx \bigsqcup_{t \in \u {\wt T}} S_t.
\end{equation}
However, if $\u T=\si{a}{c}$ and $\u {\wt T}=\si{a}{b} \sqcup \si{b+1}{c}$ (see Problem~\ref{prob:alpha}), then the identity is true because either $\si{a}{b}$ and
$\si{b+1}{c}$ are disjoint (in which case the two boxes have the same sign) or one set is contained in the other (in which case the two boxes have opposite sign).}

\section{Some sijections on signed boxes} \label{sec:sb}

The first sijection in this section will serve as the base of induction for Problem \ref{prob:pi}.

\comment{

\begin{example} \label{ex:toempty}
  For $a,b \in \Z$, define $\u T = \si{a+1}{b+1} \times \si a b$. Then $\psi: \u T \Rightarrow \u T$, defined by $\psi((l_1,l_2),i) = ((l_2+1,l_1-1),1-i)$ for $l_1 \in \si{a+1}{b+1}, l_2 \in \si a b, i \in \{0,1\}$, is a (non-normal) sijection. Define $\u S_{((l_1,l_2),i)} = (-1)^i \si{l_1}{l_2}$. Since $\u S_{((l_2+1,l_1-1),i)} = - \u S_{((l_1,l_2),i)}$ and $\u S_{((l_1,l_2),1-i)} = - \u S_{((l_1,l_2),i)}$, we have $\u S_{\psi((l_1,l_2),i)} = \u S_{((l_1,l_2),i)}$. Therefore we can use Proposition \ref{prop:sijections} with $\varphi_t = \id$, and we obtain a sijection
  $$\bigsqcup_{(l_1,l_2) \in \si{a+1}{b+1} \times \si a b} \si{l_1}{l_2} \Rightarrow \bigsqcup_{(l_1,l_2) \in \si{a+1}{b+1} \times \si a b} -\si{l_1}{l_2},$$
  which is equivalent to a sijection
  $$\bigsqcup_{(l_1,l_2) \in \si{a+1}{b+1} \times \si a b} \si{l_1}{l_2} \Rightarrow \u \emptyset.$$
\end{example}}

\begin{example} \label{ex:toempty}
  For $a,b \in \Z$, we have a normal sijection
  $$\bigsqcup_{(l_1,l_2) \in \si{a+1}{b+1} \times \si a b} \si{l_1}{l_2} \Rightarrow \u \emptyset$$
  defined by $\varphi((x,(l_1,l_2)),0) = ((x,(l_2+1,l_1-1)),0)$. It is well defined because $(l_1,l_2) \in \si{a+1}{b+1} \times \si a b$ if and only if $(l_2+1,l_1-1) \in \si{a+1}{b+1} \times \si a b$, and because $x \in \si{l_1}{l_2}$ if and only if $x \in \si{l_2+1}{l_1-1}$.
\end{example}

Note that the $0$ as the second coordinate in the example comes from the fact that a sijection in question is an involution on the disjoint union
$$\left(\bigsqcup_{(l_1,l_2) \in \si{a+1}{b+1} \times \si a b} \si{l_1}{l_2}\right) \sqcup \u \emptyset = \left(\bigsqcup_{(l_1,l_2) \in \si{a+1}{b+1} \times \si a b} \si{l_1}{l_2}\right) \times \u{\{0\}} \cup \u \emptyset \times \u{\{1\}}.$$
We could be a little less precise and write $\varphi(x,(l_1,l_2)) = (x,(l_2+1,l_1-1))$ without causing confusion.

\medskip

The following generalizes the construction of Problem \ref{prob:alpha}; indeed, for $n = 2$ the construction gives a sijection from $[a_1,b_1]$ to $[a_1,x] \sqcup (-[b_1+1,x])$.

\begin{problem} \label{prob:beta}
 Given $\q a =(a_1,\ldots,a_{n-1}) \in \Z^{n-1}$, $\q b=(b_1,\ldots,b_{n-1}) \in \Z^{n-1}$, $x \in \Z$, construct a normal sijection
 $$ \beta = \beta_{\q a,\q b,x} \colon \si{a_1}{b_1} \times \cdots \times \si{a_{n-1}}{b_{n-1}} \Rightarrow \bigsqcup_{\q (l_1,\ldots,l_{n-1}) \in \u S_1 \times \cdots \times \u S_{n-1}} \si{l_1}{l_2} \times \si{l_2}{l_3} \times \cdots \times \si{l_{n-2}}{l_{n-1}} \times \si{l_{n-1}} x,$$
 where $\u S_i = (\{a_i\},\emptyset) \sqcup (\emptyset,\{b_i+1\})$
\end{problem}

{Note that $(\{a_i\},\emptyset) \sqcup (\emptyset,\{b_i+1\})$ can be identified with
$(\{a_i\},\{b_i+1\})$ if $a_i \not= b_i+1$.}

\begin{proof}[Construction]
 The proof is by induction, with the case $n = 1$ being trivial and the case $n=2$ was constructed in Problem~\ref{prob:alpha}. Now, for $n \ge 3$,
 \begin{multline*}
   \si{a_1}{b_1} \times \cdots \times \si{a_{n-1}}{b_{n-1}} \approx \si{a_1}{b_1} \times \bigsqcup_{\q (l_2,\ldots,l_{n-1}) \in \u S_2 \times \cdots \times \u S_{n-1}} \si{l_2}{l_3} \times \cdots \times \si{l_{n-2}}{l_{n-1}} \times \si{l_{n-1}}x \\
   \approx \left( \si{a_1}{b_1} \times \bigsqcup_{\q (l_3,\ldots,l_{n-1}) \in \u S_3 \times \cdots \times \u S_{n-1}} \si{a_2}{l_3} \times \cdots \times \si{l_{n-1}}x \right) \\ \sqcup \left( \si{a_1}{b_1} \times \bigsqcup_{\q (l_3,\ldots,l_{n-1}) \in \u S_3 \times \cdots \times \u S_{n-1}} (-\si{b_2+1}{l_3}) \times \cdots \times \si{l_{n-1}}x \right),
\end{multline*}
where we used induction for the first equivalence, and distributivity and the fact that $S_2 = (\{a_2\},\emptyset) \sqcup (\emptyset,\{b_2+1\})$ for the second equivalence. By Problem~\ref{prob:alpha} and Proposition~\ref{prop:sijections} (2), there exists a sijection from the last expression to
\begin{multline*}
   \left( \left(\si{a_1}{a_2} \sqcup (-\si{b_1+1}{a_2}) \right) \times \bigsqcup_{\q (l_3,\ldots,l_{n-1}) \in \u S_3 \times \cdots \times \u S_{n-1}} \si{a_2}{l_3} \times \cdots \times \si{l_{n-1}}x \right) \\
   \sqcup \left( \left( \si{a_1}{b_2+1} \sqcup (-\si{b_1+1}{b_2+1}) \right) \times \bigsqcup_{\q (l_3,\ldots,l_{n-1}) \in \u S_3 \times \cdots \times \u S_{n-1}} (-\si{b_2+1}{l_3}) \times \cdots \times \si{l_{n-1}}x \right) \\
   \approx \bigsqcup_{\q (l_1,\ldots,l_{n-1}) \in \u S_1 \times \cdots \times \u S_{n-1}} \si{l_1}{l_2} \times \si{l_2}{l_3} \times \cdots \cdots \si{l_{n-2}}{l_{n-1}} \times \si{l_{n-1}} x,
\end{multline*}
where for the last equivalence we have again used distributivity. Normality follows from the normality of all the sijections involved in the construction.
\end{proof}

\begin{problem} \label{prob:gamma}
 Given $\q k =(k_1,\ldots,k_n) \in \Z^n$ and $x \in \Z$, construct a normal sijection
 \begin{multline*}
   \gamma = \gamma_{\q k,x} \colon \si{k_1}{k_2} \times \cdots \times \si{k_{n-1}}{k_n} \\
   \Rightarrow \bigsqcup_{i = 1}^n \si{k_1}{k_2} \times \cdots \times \si{k_{i-1}}{x+n-i} \times \si{x+n-i}{k_{i+1}} \times \cdots \times \si{k_{n-1}}{k_n} \\
   \sqcup \bigsqcup_{i = 1}^{n-2} \cdots \times \si{k_{i-1}}{k_i} \times \si{k_{i+1}+1}{x+n-i-1} \times \si{k_{i+1}}{x+n-i-2} \times \si{k_{i+2}}{k_{i+3}} \times\cdots .
   \end{multline*}
 \end{problem}

\begin{proof}[Construction] The proof is by induction with respect to $n$. The case $n=1$ is trivial, and $n=2$ is Problem~\ref{prob:alpha}. Now take $n > 2$. By the induction hypothesis (for $(k_1,\ldots,k_{n-1})$ and $x+1$), we have
\begin{multline*}\si{k_1}{k_2} \times \cdots \times \si{k_{n-1}}{k_n} \approx \bigg( \bigsqcup_{i = 1}^{n-1} \si{k_1}{k_2} \times \cdots \times \si{k_{i-1}}{x+n-i} \times \si{x+n-i}{k_{i+1}} \times \cdots \times \si{k_{n-2}}{k_{n-1}}\\
  \sqcup \bigsqcup_{i = 1}^{n-3} \si{k_1}{k_2} \times \cdots \times \si{k_{i+1}+1}{x+n-i-1} \times \si{k_{i+1}}{x+n-i-2} \times \cdots \times \si{k_{n-2}}{k_{n-1}} \bigg) \times \si{k_{n-1}}{k_n}.
  \end{multline*}
We use distributivity. We keep all terms except the one corresponding to $i = n-1$ in the first part. Because
\begin{multline*}
\si{k_{n-2}}{x+1} \times \si{k_{n-1}}{k_n} \approx \si{k_{n-2}}{x+1} \times (\si{k_{n-1}}x \sqcup \si{x+1}{k_n}) \\
\approx (\si{k_{n-2}}{k_{n-1}} \sqcup \si{k_{n-1}+1}{x+1}) \times \si{k_{n-1}}x \sqcup \si{k_{n-2}}{x+1} \times \si{x+1}{k_n} \\
{\approx \si{k_{n-2}}{k_{n-1}} \times \si{k_{n-1}}x \sqcup \si{k_{n-1}+1}{x+1} \times \si{k_{n-1}}x \sqcup \si{k_{n-2}}{x+1} \times \si{x+1}{k_n}},
\end{multline*}
we obtain the required Cartesian products for the first term on the right-hand side at $i = n$, the second term at $i = n-2$, and the first term at $i = n-1$. Again, normality follows from the fact that $\alpha$ is normal.
\end{proof}

\section{Gelfand-Tsetlin patterns} \label{sec:gt}

Using our definition of a disjoint union of {signed sets}, it is easy to define generalized Gelfand-Tsetlin patterns, or GT patterns for short (compare with \cite{Fis05}).

\begin{definition}
 For $k \in \Z$, define ${\GT(k) = (\{\cdot\},\emptyset)}$,  
 and for $\q k = (k_1,\ldots,k_n) \in \Z^n$, define {recursively}
 $$\GT(\q k) = \GT(k_1,\ldots,k_n) = \bigsqcup_{\q l \in \si{k_1}{k_2} \times \cdots \times \si{k_{n-1}}{k_n}} \GT(l_1,\ldots,l_{n-1}).$$
\end{definition}

In particular, $\GT(a,b) \approx \si a b$.

\medskip

Of course, one can think of an element of $\GT(\q k)$ in the usual way, as a triangular array $A=(A_{i,j})_{1 \leq j \leq i \le n}$ of $\binom {n+1}2$ numbers, 
arranged as
$$
\begin{array}{ccccccccc}
 &&&& A_{1,1} &&&& \\
 &&& A_{2,1} && A_{2,2} &&& \\
 && A_{3,1} && A_{3,2} && A_{3,3} && \\
 & \iddots & \vdots & \ddots & \vdots & \iddots & \vdots & \ddots &  \\
A_{n,1} && A_{n,2} && \ldots && \ldots && A_{n,n},
 \end{array}
$$
so that $A_{i+1,j} \leq A_{i,j} \leq A_{i+1,j+1}$ or $A_{i+1,j} > A_{i,j} > A_{i+1,j+1}$ for $1 \leq j \leq i < n$, and $A_{n,i}=k_i$. The sign of such an array is $(-1)^m$, where $m$ is the number of $(i,j)$ with $a_{i,j} > a_{i,j+1}$.

\medskip

Some crucial sijections for GT patterns are given by the following constructions.

\begin{problem} \label{prob:rho}
  Given $\q a =(a_1,\ldots,a_{n-1}) \in \Z^{n-1}$, $\q b=(b_1,\ldots,b_{n-1}) \in \Z^{n-1}$, $x \in \Z$, construct a sijection
  $$ \rho = \rho_{\q a,\q b,x} \colon \bigsqcup_{\q l \in \si{a_1}{b_1} \times \cdots \times \si{a_{n-1}}{b_{n-1}}} \GT(\q l) \Rightarrow \bigsqcup_{\q (l_1,\ldots,l_{n-1}) \in \u S_1 \times \cdots \times \u S_{n-1}} \GT(l_1,\ldots,l_{n-1},x),$$
  where $\u S_i = (\{a_i\},\emptyset) \sqcup (\emptyset,\{b_i+1\})$.
\end{problem}
\begin{proof}[Construction]
 In Problem \ref{prob:beta}, we constructed a normal sijection
 $$ \si{a_1}{b_1} \times \cdots \times \si{a_{n-1}}{b_{n-1}} \Rightarrow \bigsqcup_{\q (l_1,\ldots,l_{n-1}) \in \u S_1 \times \cdots \times \u S_{n-1}} \si{l_1}{l_2} \times \si{l_2}{l_3} \times \cdots \times \si{l_{n-2}}{l_{n-1}} \times \si{l_{n-1}} x.$$
 By Proposition~\ref{prop:sijections} (3) (see the comment at the end of Section \ref{sec:ss}), this gives a sijection
 $$ \bigsqcup_{\q l \in \si{a_1}{b_1} \times \cdots \times \si{a_{n-1}}{b_{n-1}}} \GT(\q l) \Rightarrow \bigsqcup_{\q m \in \bigsqcup_{\q (l_1,\ldots,l_{n-1}) \in \u S_1 \times \cdots \times \u S_{n-1}} \si{l_1}{l_2} \times \si{l_2}{l_3} \times \cdots \times \si{l_{n-2}}{l_{n-1}} \times \si{l_{n-1}} x} \GT(\q m).$$
 By basic sijection constructions, we get that this is equivalent to
 $$\bigsqcup_{\q (l_1,\ldots,l_{n-1}) \in \u S_1 \times \cdots \times \u S_{n-1}} \bigsqcup_{\q m \in \si{l_1}{l_2} \times \si{l_2}{l_3} \times \cdots \times \si{l_{n-2}}{l_{n-1}} \times \si{l_{n-1}} x} \GT(\q m),$$
 and by definition of $\GT$, this is equal to $\bigsqcup_{\q (l_1,\ldots,l_{n-1}) \in \u S_1 \times \cdots \times \u S_{n-1}} \GT(l_1,\ldots,l_{n-1},x)$.
\end{proof}

The result is important because while it adds a dimension to GT patterns, it (typically) greatly reduces the size of the indexing signed set. {In fact, there is an analogy to the fundamental theorem of calculus: instead of extending the disjoint union over the entire signed box, it suffices to consider the boundary; $x$ corresponds in a sense to the constant of integration.}

\begin{problem} \label{prob:pi}
  Given $\q k =(k_1,\ldots,k_n) \in \Z^n$ and $i$, $1 \leq i \leq n-1$, construct a sijection
 $$\pi = \pi_{\q k,i} \colon \GT(k_1,\ldots,k_n) \Rightarrow -\GT(k_1,\ldots,k_{i-1},k_{i+1}+1,k_i-1,k_{i+2},\ldots,k_n).$$
 Given $\q a =(a_1,\ldots,a_n) \in \Z^n$,  $\q b =(b_1,\ldots,b_n) \in \Z^n$ such that for some $i$, $1 \leq i \leq n-1$, we have $a_{i+1} = a_i - 1$ and $b_{i+1} = b_i - 1$, construct a sijection
 $$\sigma = \sigma_{\q a,\q b,i} \colon \bigsqcup_{\q l \in \si{a_1}{b_1} \times \cdots \times \si{a_n}{b_n}} \GT(\q l) \Rightarrow \u \emptyset.$$
\end{problem}
\begin{proof}[Construction]
 The proof is by induction, with the induction step for $\pi$ using $\sigma$ and vice versa. For $n = 1$, there is nothing to prove. For $n=2$ and $i=1$, the existence of $\pi$ follows from the statement $\si{k_1}{k_2} = -\si{k_2+1}{k_1-1}$, and $\sigma$ was constructed in Example \ref{ex:toempty}. Assume that $n > 2$ and $1 < i < n-1$. We have
$$\GT(k_1,\ldots,k_n) =  \bigsqcup_{\q l \in \si{k_1}{k_2} \times \cdots \times \si{k_{i-1}}{k_i} \times \si{k_{i+1}+1}{k_i-1} \times \si{k_{i+1}}{k_{i+2}} \times \cdots \times \si{k_{n-1}}{k_n}} -\GT(l_1,\ldots,l_{n-1}).$$
By using $\id \times \cdots \times \id \times \alpha_{k_{i-1},k_{i+1}+1,k_i} \times \id \times \alpha_{k_{i+1},k_i-2,k_{i+2}} \times \id \times \cdots \times \id$ and distributivity, we get a normal sijection
\begin{multline*}
\si{k_1}{k_2} \times \cdots \times \si{k_{i-1}}{k_i} \times \si{k_{i+1}+1}{k_i-1} \times \si{k_{i+1}}{k_{i+2}} \times \cdots \times \si{k_{n-1}}{k_n} \Rightarrow \\
\si{k_1}{k_2} \times \cdots \times \si{k_{i-1}}{k_{i+1}+1} \times \si{k_{i+1}+1}{k_i-1} \times \si{k_{i}-1}{k_{i+2}} \times \cdots \times \si{k_{n-1}}{k_n} \\
\sqcup \si{k_1}{k_2} \times \cdots \times \si{k_{i-1}}{k_{i+1}+1} \times \si{k_{i+1}+1}{k_i-1} \times \si{k_{i+1}}{k_{i}-2} \times \cdots \times \si{k_{n-1}}{k_n} \\
\sqcup \si{k_1}{k_2} \times \cdots \times \si{k_{i+1}+2}{k_i} \times \si{k_{i+1}+1}{k_i-1} \times \si{k_{i}-1}{k_{i+2}} \times \cdots \times \si{k_{n-1}}{k_n} \\
\sqcup \si{k_1}{k_2} \times \cdots \times \si{k_{i+1}+2}{k_i} \times \si{k_{i+1}+1}{k_i-1} \times \si{k_{i+1}}{k_{i}-2} \times \cdots \times \si{k_{n-1}}{k_n}
\end{multline*}
 By Proposition~\ref{prop:sijections} (3), this gives a sijection
 \begin{multline*}
  \bigsqcup_{\q l \in \si{k_1}{k_2} \times \cdots \times \si{k_{i-1}}{k_i} \times \si{k_{i+1}+1}{k_i-1} \times \si{k_{i+1}}{k_{i+2}} \times \cdots \times \si{k_{n-1}}{k_n}} -\GT(l_1,\ldots,l_{n-1}) \Rightarrow \\
   \bigsqcup_{\q l \in \si{k_1}{k_2} \times \cdots \times \si{k_{i-1}}{k_{i+1}+1} \times \si{k_{i+1}+1}{k_i-1} \times \si{k_{i}-1}{k_{i+2}} \times \cdots \times \si{k_{n-1}}{k_n}} - \GT(l_1,\ldots,l_{n-1}) \\ \sqcup  \bigsqcup_{\q l \in \si{k_1}{k_2} \times \cdots \times \si{k_{i-1}}{k_{i+1}+1} \times \si{k_{i+1}+1}{k_i-1} \times \si{k_{i+1}}{k_{i}-2} \times \cdots \times \si{k_{n-1}}{k_n}} -\GT(l_1,\ldots,l_{n-1})  \\ \sqcup  \bigsqcup_{\q l \in \si{k_1}{k_2} \times \cdots \times \si{k_{i+1}+2}{k_i} \times \si{k_{i+1}+1}{k_i-1} \times \si{k_{i}-1}{k_{i+2}} \times \cdots \times \si{k_{n-1}}{k_n}} -\GT(l_1,\ldots,l_{n-1})
 \\ \sqcup  \bigsqcup_{\q l \in \si{k_1}{k_2} \times \cdots \times \si{k_{i+1}+2}{k_i} \times \si{k_{i+1}+1}{k_i-1} \times \si{k_{i+1}}{k_{i}-2} \times \cdots \times \si{k_{n-1}}{k_n}} -\GT(l_1,\ldots,l_{n-1}).
 \end{multline*}
 By definition, the first signed set on the right-hand side is $-\GT(k_1,\ldots,k_{i-1},k_{i+1}+1,k_i-1,k_{i+2},\ldots,k_n)$. The other three disjoint unions all satisfy the condition needed for the existence of $\sigma$ {(for $i$, for $i-1$ and for both, $i-1$ and $i$, respectively)}, and hence we can siject them to $\u \emptyset$.\\
 If $i = 1$ or $i = n-1$, the proof is similar but easier (as we only have to use $\alpha$ once, and we get only two factors after using distributivity). Details are left to the reader.\\
 Now take $\q l = (l_1,\ldots,l_n)$ and $\q{l'}=(l_1,\ldots,l_{i-1},l_{i+1}+1,l_i-1,l_{i+2},\ldots,l_n)$. The sijection $\sigma$ can then be defined as
 $$\sigma_{\q a,\q b,i}(A,\q l) = \begin{cases} (\pi_{\q l,i}(A),\q l) & \mbox{if }\pi_{\q l,i}(A) \in \GT(\q l) \\ (\pi_{\q l,i}(A),\q {l'}) & \mbox{if } \pi_{\q {l},i}(A) \in \GT(\q {l'})  \end{cases}.$$
 It is easy to check that this is a sijection. Compare with Example \ref{ex:toempty}.
\end{proof}

\begin{problem} \label{prob:tau}
 Given $\q k =(k_1,\ldots,k_n) \in \Z^n$ and $x \in \Z$, construct a sijection
 $$\tau = \tau_{\q k,x} \colon \GT(k_1,\ldots,k_n) \Rightarrow \bigsqcup_{i = 1}^n \GT(k_1,\ldots,k_{i-1},x+n-i,k_{i+1},\ldots,k_n).$$
\end{problem}
\begin{proof}[Construction]
 In Problem \ref{prob:gamma}, we constructed a normal sijection
 \begin{multline*}
   \si{k_1}{k_2} \times \cdots \times \si{k_{n-1}}{k_n} \Rightarrow \bigsqcup_{i = 1}^n \si{k_1}{k_2} \times \cdots \times \si{k_{i-1}}{x+n-i} \times \si{x+n-i}{k_{i+1}} \times \cdots \times \si{k_{n-1}}{k_n}\\
   \sqcup \bigsqcup_{i = 1}^{n-2} \cdots \times \si{k_{i-1}}{k_i} \times \si{k_{i+1}+1}{x+n-i-1} \times \si{k_{i+1}}{x+n-i-2} \times \si{k_{i+2}}{k_{i+3}} \times\cdots ,
   \end{multline*}
which gives a sijection
\begin{multline*}\bigsqcup_{\q l \in \si{k_1}{k_2} \times \cdots \times \si{k_{n-1}}{k_n}} \GT(\q l) \Rightarrow \bigsqcup_{i = 1}^n \bigsqcup_{\q l \in \si{k_1}{k_2} \times \cdots \times \si{k_{i-1}}{x+n-i} \times \si{x+n-i}{k_{i+1}} \times \cdots \times \si{k_{n-1}}{k_n}}\GT(\q l) \\
\sqcup \bigsqcup_{i = 1}^{n-2} \bigsqcup_{\q l \in \si{k_1}{k_2} \times \cdots \times \si{k_{i-1}}{k_i} \times \si{k_{i+1}+1}{x+n-i-1} \times \si{k_{i+1}}{x+n-i-2} \times \si{k_{i+2}}{k_{i+3}} \times\cdots \times \si{k_{n-1}}{k_n}} \GT(\q l).
\end{multline*}
All disjoint unions in the second term satisfy the conditions for the existence of $\sigma$ from Problem \ref{prob:pi}, so we can siject them to $\u \emptyset$. This gives a sijection
 $$\bigsqcup_{\q l \in \si{k_1}{k_2} \times \cdots \times \si{k_{n-1}}{k_n}} \GT(\q l) \Rightarrow \bigsqcup_{i = 1}^n \bigsqcup_{\q l \in \si{k_1}{k_2} \times \cdots \times \si{k_{i-1}}{x+n-i} \times \si{x+n-i}{k_{i+1}} \times \cdots \times \si{k_{n-1}}{k_n}}\GT(\q l),$$
 which is, by the definition of $\GT$, a sijection $\GT(k_1,\ldots,k_n) \Rightarrow \bigsqcup_{i = 1}^n \GT(k_1,\ldots,k_{i-1},x+n-i,k_{i+1},\ldots,k_n)$.
 \comment{Finally, for each $i$, apply the sijection
 $$\pi_{(k_1,\ldots,k_{i-1},k_{i+1}+1,\ldots,k_{n-1}+1,x+1,k_n),n-1} \circ \cdots \circ \pi_{(k_1,\ldots,k_{i-1},k_{i+1}+1,x+n-i-1,k_{i+2},\ldots,k_n),i+1} \circ \pi_{(k_1,\ldots,k_{i-1},x+n-i,k_{i+1},\ldots,k_n),i}$$
 from $\GT(k_1,\ldots,k_{i-1},x+n-i,k_{i+1},\ldots,k_n)$ to $(-1)^{n-i} \GT(k_1,\ldots,k_{i-1},k_{i+1}+1,\ldots,k_n+1,x)$.}
\end{proof}

\section{Combinatorics of the monotone triangle recursion} \label{sec:recursion}

\subsection*{Monotone triangles}

Suppose that $\q k = (k_1,\ldots,k_n)$ and $\q l = (l_1,\ldots,l_{n-1})$ are two sequences of integers. We say that $\q l$ \emph{interlaces} $\q k$, $\q l \prec \q k$, if the following holds:
\begin{enumerate}
  \item for every $i$, $1 \leq i \leq n-1$, $l_i$ is in the closed interval between $k_i$ and $k_{i+1}$;
  \item if $k_{i-1} \leq k_i \leq k_{i+1}$ for some $i$, $2 \leq i \leq n-1$, then $l_{i-1}$ and $l_i$ cannot both be $k_i$;
  \item if $k_i > l_i = k_{i+1}$, then $i \leq n-2$ and $l_{i+1} = l_i = k_{i+1}$;
  \item if $k_i = l_i > k_{i+1}$, then $i \geq 2$ and $l_{i-1} = l_i = k_i$.
\end{enumerate}

For example, if $k_1 < k_2 < \ldots < k_n$, then $l_i \in [k_i,k_{i+1}]$ and $l_1 < l_2 < \ldots < l_{n-1}$.

\medskip

A \emph{monotone triangle of size $n$} is a map $T \colon \{(i,j) \colon 1 \leq j \leq i \leq n \} \to \Z$ so that line $i-1$ (i.e.~the sequence $T_{i-1,1},\ldots,T_{i-1,i-1}$) interlaces line $i$ (i.e.~the sequence $T_{i,1},\ldots,T_{i,i}$).
\begin{example} {The following is a monotone triangles of size $5$:}
$$\begin{array}{ccccccccc}
&&&& 4 &&&& \\
&&& 3 && 5 &&& \\
&& 3 && 4 && 5 && \\
& 3 & & 3 & & 4 & & 5 & \\
5 & & 3 & & 1 & & 4 & & 6
\end{array}
$$
\end{example}
{This notion of (generalized) monotone triangle was introduced in \cite{Rie13}. Other notions appeared in \cite{Fis12}.} 

\medskip

The \emph{sign} of a monotone triangle $T$ is $(-1)^r$, where $r$ is the sum of:
\begin{itemize}
  \item the number of strict descents in the rows of $T$, i.e.~the number of pairs $(i,j)$ so that $1 \leq j < i \leq n$ and $T_{i,j} > T_{i,j+1}$, and
  \item the number of $(i,j)$ so that $1 \leq j \leq i - 2$, $i \leq n$ and $T_{i,j} > T_{i-1,j} = T_{i,j+1} = T_{i-1,j+1} > T_{i,j+2}$.
\end{itemize}
{The sign of our example is $-1$.}

\medskip

We denote the signed set of all monotone triangles with bottom row $\q k$ by $\mtk$.

\medskip

It turns out that $\MT(\q k)$ satisfies a recursive ``identity''. Let us define the signed set of \emph{arrow rows of order $n$} as
$$\AR_n = (\{\nearrow,\nwarrow\},\{\nenwarrow\})^{n}.$$
Alternatively, we can think of them as rows of length $n$ with elements $\nwarrow, \nearrow, \nenwarrow$, where the positive elements are precisely those with an even number of $\nenwarrow$'s.

\medskip

The role of an arrow row $\mu$ of order $n$ is that it induces a deformation of $\si{k_1}{k_2} \times \si{k_2}{k_3} \times \cdots \times \si{k_{n-1}}{k_n}$ as follows. Consider
$$
\begin{array}{ccccccccccc}
 & \si{k_1}{k_2} &   & \si{k_2}{k_3} &  &  \ldots &  & \si{k_{n-2}}{k_n-1} &  & \si{k_{n-1}}{k_n} & \\
 \mu_1 & & \mu_2  & & \mu_3 &  &  \ldots &   & \mu_{n-1} & & \mu_{n} ,
 \end{array}
 $$
and if $\mu_i \in \{\nwarrow, \nenwarrow\}$ (that is we have an arrow pointing towards $\si{k_{i-1}}{k_i}$) then
$k_i$ is decreased by $1$ in $\si{k_{i-1}}{k_i}$, while there is no change for this $k_i$ if $\mu_i=\nearrow$. If $\mu_i \in \{\nearrow, \nenwarrow\}$ (that is we have an arrow pointing towards $[k_i,k_{i+1}]$) then $k_i$ is increased by $1$ in $[k_i,k_{i+1}]$, while there is no change for this $k_i$ if $\mu_i=\nwarrow$.

\medskip

For a more formal description, we let $\delta_{\nwarrow}(\nwarrow) = \delta_{\nwarrow}(\nenwarrow) = \delta_{\nearrow}(\nearrow) = \delta_{\nearrow}(\nenwarrow) = 1$ and $\delta_{\nwarrow}(\nearrow) = \delta_{\nearrow}(\nwarrow) = 0$, and we define
$$e(\q k,\mu)  = \si{k_1+\delta_{\nearrow}(\mu_1)}{k_2-\delta_{\nwarrow}(\mu_2)} \times \ldots \times \si{k_{n-1}+\delta_{\nearrow}(\mu_{n-1})}{k_n-\delta_{\nwarrow}(\mu_n)}.$$
for $\q k = (k_1,\ldots,k_n)$ and $\mu \in \AR_n$.

\begin{problem} \label{prob:Xi}
  Given $\q k = (k_1,\ldots,k_n)$, construct a sijection
  $$\Xi = \Xi_{\q k} \colon \MT(\q k) \Rightarrow \bigsqcup_{\mu \in \AR_n} \bigsqcup_{\q l \in e(\q k,\mu)} \MT(\q l).$$
\end{problem}
\begin{proof}[Construction]
 All elements on the left are mapped to the right with $\Xi$, while there are quite a few cancellations on the right-hand side. More specifically, take a monotone triangle $T$ with bottom row $\q k$. Then $\Xi(T) = ((T',\q l),\mu)$, where $T'$ is the monotone triangle we obtain from $T$ by deleting the last row, $\q l$ is the bottom row of $T'$, and $\mu = (\mu_1,\ldots,\mu_n)$ is the arrow row defined as follows:
 \begin{itemize}
  \item $\mu_1 = \nwarrow$;
  \item $\mu_n = \nearrow$;
  \item for $ 1 < i < n$, $\mu_i$ is determined as follows:
  \begin{enumerate} \item if $k_{i-1} \leq l_{i-1} = k_i$, take $\mu_i = \nearrow$;  \item if $k_{i-1} > l_{i-1} = k_i = l_{i} > k_{i+1}$, take $\mu_i = \nenwarrow$; \item otherwise, take $\mu_i = \nwarrow$. \end{enumerate}
 \end{itemize}
 It is easy to check that $\q l$ is indeed in $e(\q k,\mu)$.
 Note that in (1) and (2) of the third bullet point, $\mu_i$ is forced if we require $\q l \in e(\q k,\mu)$. In (3), $\mu_i=\nenwarrow$ would also be possible if and only if  $\mu_i=\nearrow$ would also be possible.\\
 On the other hand, for $((T',\q l),\mu)$, define $\Xi((T',\q l),\mu)$ as follows. For the construction it is useful to keep in mind that $\q l \in e(\q k, \mu)$ implies that conditions (1) and (2) for $\q l \prec \q k$ are satisfied.
 \begin{itemize}
  \item if $\mu_1 \neq \nwarrow$, take $\Xi((T',\q l),\mu) = ((T',\q l),\mu')$, where we obtain $\mu'$ from $\mu$ by replacing $\nenwarrow$ in position $1$ by $\nearrow$ and vice versa;
  \item if $\mu_1 = \nwarrow$ and $\mu_n \neq \nearrow$, take $\Xi((T',\q l),\mu) = ((T',\q l),\mu')$, where we obtain $\mu'$ from $\mu$ by replacing $\nenwarrow$ in position $n$ by $\nwarrow$ and vice versa;
  \item if $\mu_1 = \nwarrow$ and $\mu_n = \nearrow$, and $\q l \not\prec \q k$, find the smallest $i$ between $2$ and $n-1$ such that:
 \begin{itemize}
   \item condition (3) of $\q l \prec \q k$ is not satisfied at $i$, i.e.
   $k_{i-1} > l_{i-1} = k_{i} \not= l_{i}$ (which implies $\mu_i \in \{\nwarrow, \nenwarrow\}$), or
   \item condition (4) of $\q l \prec \q k$ is not satisfied at $i$, i.e.
   $l_{i-1} \not= k_i = l_i > k_{i+1}$ (which implies  $\mu_i \in \{\nearrow, \nenwarrow\}$).
\end{itemize}
Then take $\Xi((T',\q l),\mu) = ((T',\q l),\mu')$, where we obtain $\mu'$ from $\mu$ by replacing $\nenwarrow$ in position $i$ by $\nwarrow$ and vice versa in the first case, and replacing $\nenwarrow$ in position $i$ by $\nearrow$ and vice versa in the second case;
 \item if $\mu_1 = \nwarrow$ and $\mu_n = \nearrow$, and $\q l \prec \q k$, find the smallest $i$ for an instance of (3) of the third bullet point in the first paragraph of the proof with $\mu_i \not= \nwarrow$ (if such an $i$ exists). Then take $\Xi((T',\q l),\mu) = ((T',\q l),\mu')$, where we obtain $\mu'$ from $\mu$ by replacing $\nenwarrow$ in position $i$ by $\nwarrow$ and vice versa.
 \end{itemize}
If no such $i$ exists, we take $\Xi((T',\q l),\mu) = T$, where we obtain $T$ from $T'$ by adding $\q k$ as the last row. It is easy to see that this is a well-defined sijection.
\end{proof}

\begin{remark} The previous construction could have been avoided by using alternative extensions of monotone triangles provided in \cite{Fis12}. However, the advantage of the definition used in this paper is that it is more reduced than the others in the sense that it can obtained from these by cancelling elements using certain sign-reversing involutions.
\end{remark}

\subsection*{Arrow patterns and shifted GT patterns}

Define the signed set of \emph{arrow patterns of order $n$} as
$$\AP_n = (\{\swarrow, \searrow\},\{\seswarrow\})^{\binom n 2}.$$

\medskip

Alternatively, we can think of an arrow pattern of order $n$ as a triangular array $T=(t_{p,q})_{1 \le p < q \le n}$ arranged as
$$T = \begin{smallmatrix} & & & & t_{1,n} & & & & \\ & & & t_{1,n-1} & & t_{2,n} & & & \\ & & t_{1,n-2} & & t_{2,n-1}& & t_{3,n} & & \\ & \vstretch{0.35}{\udots} & \vstretch{0.35}{\vdots} & \vstretch{0.35}{\ddots} & \vstretch{0.35}{\vdots} & \vstretch{0.35}{\udots} & \vstretch{0.35}{\vdots} & \vstretch{0.35}{\ddots} & \\ t_{1,2} & & t_{2,3}  & & \ldots  & & \ldots  & & t_{n-1,n}  \end{smallmatrix},$$
with $t_{p,q} \in \{\swarrow, \searrow, \seswarrow\}$, and the sign of an arrow pattern is $1$ if the number of $\seswarrow$'s is even and $-1$ otherwise.

\medskip

The role of an arrow pattern of order $n$ is that it induces a deformation of $(k_1,\ldots,k_n)$, which can be thought of as follows. Add $k_1,\ldots,k_n$ as bottom row of $T$ (i.e., $t_{i,i}=k_i$), and for each $\swarrow$ or$ \seswarrow$ which is in the same $\swarrow$-diagonal as $k_i$ add $1$ to $k_i$, while for each $\searrow$ or $\seswarrow$ which is in the same $\searrow$-diagonal as $k_i$ subtract $1$ from $k_i$.
More formally, letting $\delta_{\swarrow}(\swarrow) = \delta_{\swarrow}(\seswarrow) = \delta_{\searrow}(\searrow) = \delta_{\searrow}(\seswarrow) = 1$ and $\delta_{\swarrow}(\searrow) = \delta_{\searrow}(\swarrow) = 0$, we set
$$c_i(T) = \sum_{j=i+1}^{n} \delta_{\swarrow}(t_{i,j}) - \sum_{j=1}^{i-1} \delta_{\searrow}(t_{j,i}) \mbox{ and } d(\q k,T) = (k_1+c_1(T),k_2+c_2(T),\ldots,k_n+c_n(T))$$
for $\q k = (k_1,\ldots,k_n)$ and $T \in \AP_n$.

\medskip

For $\q k =(k_1,\ldots,k_n)$ define \emph{shifted Gelfand-Tsetlin patterns}, or SGT patterns for short, as the following disjoint union of GT patterns over arrow patterns of order $n$:
$$
\SGT(\q k) = \bigsqcup_{T \in \AP_n} \GT(d(\q k,T))
$$


\medskip

Considering that $|(\{\swarrow, \searrow\},\{\seswarrow\})| = 1$ and therefore $|\AP_n| = 1$, the following is not surprising.

\begin{problem} \label{prob:Psi}
 Given $n$ and $i$, $1 \leq i \leq n$, construct a sijection
 $$\Psi = \Psi_{n,i} \colon  \AP_{n-1} \Rightarrow \AP_n.$$
\end{problem}
\begin{proof}[Construction]
 For $T \in \AP_{n-1}$, take $\Psi(T) = (t'_{p,q})_{1 \leq p < q \leq n}$ to be the arrow pattern defined by
 $$t'_{p,q} = \begin{cases} t_{p,q} & \mbox{if } p < q < i \\ t_{p,q-1} & \mbox{if } p < i < q \\ t_{p-1,q-1} & \mbox{if } i < p < q \\ \searrow & \mbox{if } p < q = i \\ \swarrow & \mbox{if } i = p < q \end{cases}.$$
{An example for $n=6$ and $i=4$ is
$$
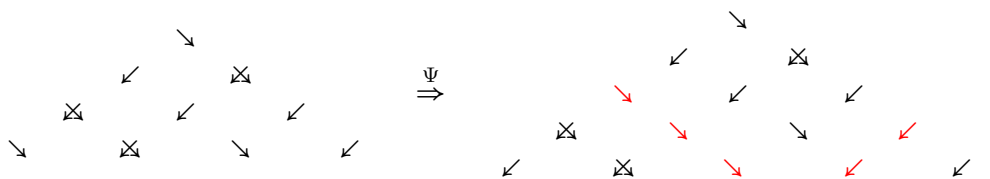
\begin{array}{ccccccc} 
&&& \searrow &&& \\
&& \swarrow && \seswarrow && \\
& \seswarrow && \swarrow && \swarrow & \\
\searrow && \seswarrow && \searrow && \swarrow
\end{array} \quad \stackrel{\Psi}{\Rightarrow} \quad
\begin{array}{ccccccccc} 
&&&& \searrow &&&& \\
&&& \swarrow && \seswarrow &&& \\
&& {\color{red} \searrow} && \swarrow &&  \swarrow && \\
& \seswarrow &&  {\color{red} \searrow} && \searrow && {\color{red} \swarrow}  & \\
\swarrow && \seswarrow && {\color{red} \searrow} && {\color{red} \swarrow} && \swarrow
\end{array},
$$
where the new arrows are indicated in red.}
 If $T \in \AP_{n}$, $t_{p,i} = \searrow$ for $p = 1,\ldots,i-1$, $t_{i,q} = \swarrow$ for $q = i+1,\ldots,n$, take $\Psi(T) = (t'_{p,q})_{1 \leq p < q \leq n-1}$, where
 $$t'_{p,q} = \begin{cases} t_{p,q} & \mbox{if } p < q < i \\ t_{p,q+1} & \mbox{if } p < i \leq q \\ t_{p+1,q+1} & \mbox{if } i \leq p < q \end{cases}.$$
Otherwise, there either exists $p$ so that $t_{p,i} \neq \searrow$, or there exists $q$ so that $t_{i,q} \neq \swarrow$. In the first case, define $\Psi(T) = (t'_{p,q})_{1 \leq p < q \leq n}$, where $t'_{p,i}=\swarrow$ if $t_{p,i} = \seswarrow$ and $t'_{p,i}=\seswarrow$ if $t_{p,i} = \swarrow$, and all other array elements are equal. In the second case, define $\Psi(T) = (t'_{p,q})_{1 \leq p < q \leq n}$, where $t'_{i,q}=\searrow$ if $t_{i,q} = \seswarrow$ and $t'_{i,q}=\seswarrow$ if $t_{i,q} = \searrow$, and all other array elements are equal. It is easy to see that this is a sijection.
\end{proof}

The difficult part of this paper is to prove that $\SGT$ satisfies the same ``recursion'' as $\MT$. While the proof of the recursion was easy for monotone triangles, it is very involved for shifted GT patterns, and needs almost {all} the sijections we have constructed in this and previous sections.

\begin{problem}
Given $\q k = (k_1,\ldots,k_n) \in \Z^n$ and $x \in \Z$, construct a sijection
$$\Phi = \Phi_{\q k,x} \colon \bigsqcup_{\mu \in \AR_n} \bigsqcup_{\q l \in e(\q k,\mu)} \SGT(\q l) \Rightarrow \SGT(\q k).$$
\end{problem}
\begin{proof}[Construction]
To make the construction of $\Phi$ a little easier, we will define it as the composition of several sijections. The first one will reduce the indexing sets (from a signed box to its ``corners'') using Problem \ref{prob:rho}. The second one increases the order of the arrow patterns using the sijection from Problem \ref{prob:Psi}. The third one further reduces the indexing set (from a signed set with $2^{n-1}$ elements to $\si{1}{n}$). The last one gets rid of the arrow row and then uses Problem \ref{prob:tau}.\\
For $\mu \in \AR_n$, define $\u S_i = (\{k_i+\delta_{\nearrow}(\mu_i)\},\emptyset) \sqcup (\emptyset,\{k_{i+1}-\delta_{\nwarrow}(\mu_{i+1})+1\})$. Then $\Phi$ is the composition of sijections
\begin{multline*}
  \bigsqcup_{\mu \in \AR_n} \bigsqcup_{\q l \in e(\q k,\mu)} \SGT(\q l) \\\stackrel{\Phi_1}\Longrightarrow \bigsqcup_{\mu \in \AR_n} \bigsqcup_{T \in \AP_{n-1}} \bigsqcup_{\q m \in \u S_1 \times \cdots \times \u S_{n-1}} \GT(m_1+c_1(T),\ldots,m_{n-1}+c_{n-1}(T),x) \\
\stackrel{\Phi_2}\Longrightarrow \bigsqcup_{\mu \in \AR_n} \bigsqcup_{T \in \AP_{n}} \bigsqcup_{\q m \in \u S_1 \times \cdots \times \u S_{n-1}} \GT(m_1+c_1(T),\ldots,m_{n-1}+c_{n-1}(T),x) \\ \stackrel{\Phi_3} \Longrightarrow
\bigsqcup_{\mu \in \AR_n} \bigsqcup_{T \in \AP_{n}} \bigsqcup_{i=1}^n \GT(\ldots,k_{i-1}+\delta_{\nearrow}(\mu_{i-1})+c_{i-1}(T), x + n-i,k_{i+1}-\delta_{\nwarrow}(\mu_{i+1})+c_{i}(T),\ldots) \\
\stackrel{\Phi_4}\Longrightarrow \SGT(\q k), \hspace*{12cm}
\end{multline*}
where $\Phi_1$, $\Phi_2$, $\Phi_3$, and $\Phi_4$ are constructed as follows.\\
\emph{Construction of $\Phi_1$.} By definition of $\SGT$, we have
$$\bigsqcup_{\mu \in \AR_n} \bigsqcup_{\q l \in e(\q k,\mu)} \SGT(\q l) = \bigsqcup_{\mu \in \AR_n} \bigsqcup_{\q l \in e(\q k,\mu)} \bigsqcup_{T \in \AP_{n-1}} \GT(d(\q l,T)).$$
By switching the inner disjoint unions, we get a sijection to
$$\bigsqcup_{\mu \in \AR_n} \bigsqcup_{T \in \AP_{n-1}} \bigsqcup_{\q l \in e(\q k,\mu)}  \GT(d(\q l,T)).$$
There is an obvious sijection from this signed set to
$$\bigsqcup_{\mu \in \AR_n} \bigsqcup_{T \in \AP_{n-1}} \bigsqcup_{\q l \in d(e(\q k,\mu),T)}  \GT(\q l),$$
{by abuse of notation setting}
$${ d(\si{x_1}{y_1} \times \ldots \si{x_{n-1}}{y_{n-1}},T)=\si{x_1+c_1(T)}{y_1+c_1(T)} \times \cdots \times \si{x_{n-1}+c_{n-1}(T)}{y_{n-1}+c_{n-1}(T)}.}$$
Now for each $\mu$ and $T$, use the map $\rho$ from Problem \ref{prob:rho} for $a_i = k_i+\delta_{\nearrow}(\mu_i)+c_i(T)$, $b_i = k_{i+1}-\delta_{\nwarrow}(\mu_{i+1})+c_i(T)$, and $x$. We get a sijection to
$$\bigsqcup_{\mu \in \AR_n} \bigsqcup_{T \in \AP_{n-1}} \bigsqcup_{\q m \in \u {S}'_1 \times \cdots \u{S}'_{n-1}} \GT(m_1,\ldots,m_{n-1},x),$$
where $\u{S}'_i = (\{k_i+\delta_{\nearrow}(\mu_i)+c_i(T)\},\emptyset) \sqcup (\emptyset,\{k_{i+1}-\delta_{\nwarrow}(\mu_{i+1})+c_i(T)+1\})$. Finally, there is an obvious sijection from this signed set to
$$\bigsqcup_{\mu \in \AR_n} \bigsqcup_{T \in \AP_{n-1}} \bigsqcup_{\q m \in \u {S}_1 \times \cdots \u{S}_{n-1}} \GT(m_1+c_1(T),\ldots,m_{n-1}+c_{n-1}(T),x).$$
\emph{Construction of $\Phi_2$.} In Problem \ref{prob:Psi}, we constructed sijections $\Psi_{n,i} \colon \AP_{n-1} \Rightarrow \AP_n$. We construct $\Phi_2$ by using Proposition \ref{prop:sijections}~ (3) for $\psi = \Psi_{n,n}$, $\u T = \AP_{n-1}$, $\u{\wt T} = \AP_n$,
$$\u S_T =  \bigsqcup_{\q m \in \u {S}_1 \times \cdots \u{S}_{n-1}} \GT(m_1+c_1(T),\ldots,m_{n-1}+c_{n-1}(T),x) \quad \mbox{for } T \in \AP_{n-1} \sqcup \AP_n$$
and $\varphi_T = \u\id$. This is well defined because $c_i(T) = c_i(\Psi_{n,n}(T))$ for $T \in \AP_{n-1} \sqcup \AP_n$ and $i=1,\ldots,n-1$.\\
\emph{Construction of $\Phi_3$.} Let $\eta$ be the involution that maps $\swarrow \leftrightarrow \nearrow, \searrow \leftrightarrow \nwarrow, \seswarrow \leftrightarrow \nenwarrow$. The elements of the signed set $\u S = \u {S}_1 \times \cdots  \times \u{S}_{n-1}$ are $(n-1)$-tuples of elements that are either $(k_i+\delta_{\nearrow}(\mu_i),0)$ or $(k_{i+1}-\delta_{\nwarrow}(\mu_{i+1})+1,1)$. Define $\u S'$ as the subset of $\u S$ containing tuples of the form $(\ldots,(m_i,1),(m_{i+1},0),\ldots)$, i.e.~the ones where we choose $k_{i+1}-\delta_{\nwarrow}(\mu_{i+1})+1$ in position $i$ and $k_{i+1}+\delta_{\nearrow}(\mu_{i+1})$ in position $i+1$ for some $i$. Then we can define a sijection
$$ \bigsqcup_{\mu \in \AR_n} \bigsqcup_{T \in \AP_{n}} \bigsqcup_{\q m \in \u S'} \GT(m_1+c_1(T),\ldots,m_{n-1}+c_{n-1}(T),x) \Rightarrow \u \emptyset$$
as follows: given $\mu \in \AR_n$, $T \in \AP_n$, $\q m = (\ldots,k_{i+1}-\delta_{\nwarrow}(\mu_{i+1})+1,k_{i+1}+\delta_{\nearrow}(\mu_{i+1}),\ldots)$ (and $i$ is the smallest index where this happens), $A \in \GT(m_1+c_1(T),\ldots,m_{n-1}+c_{n-1}(T),x)$, map $(((A,\q m),T),\mu)$ to $(((A',\q m'),T'),\mu')$, where:
\begin{itemize}
 \item $A' = \pi_{i,n}(A)$;
 \item $T'$ is $T$ if $A'$ has the same bottom row as $A$; otherwise, $T'$ is obtained from $T$ by interchanging {$t_{i,j}$ and $t_{i+1,j}$ for $j>i+1$ as well as $t_{j,i}$ and $t_{j,i+1}$ for $j<i$, and setting $t'_{i,i+1} = \eta(\mu_{i+1})$;}
 \item $\mu'$ is $\mu$ if $A'$ has the same bottom row as $A$; otherwise, $\mu'$ is obtained from $\mu$ by replacing {$\mu_{i+1}$ with $\eta(t_{i,i+1})$};
 \item $\q m' = (\ldots,k_{i+1}-\delta_{\nwarrow}(\mu'_{i+1})+1,k_{i+1}+\delta_{\nearrow}(\mu'_{i+1}),\ldots)$.
\end{itemize}
{What remains is
$$
\bigsqcup_{\mu \in \AR_n} \bigsqcup_{T \in \AP_{n}} \bigsqcup_{i=1}^n (-1)^{n-i} \GT(\ldots,k_{i-1}+\delta_{\nearrow}(\mu_{i-1})+c_{i-1}(T),k_{i+1}-\delta_{\nwarrow}(\mu_{i+1})+c_{i}(T)+1,\ldots,x),
$$
and we can now apply $\pi_{i,n} \circ \pi_{i,n} \circ \dots \circ \pi_{n-1,n}$ to obtain what is claimed.}

\emph{Construction of $\Phi_4$.} By switching the order in which we do disjoint unions, we arrive at
$$\bigsqcup_{T \in \AP_{n}}  \bigsqcup_{i=1}^n \bigsqcup_{\mu \in \AR_n} \GT(\ldots,k_{i-1}+\delta_{\nearrow}(\mu_{i-1})+c_{i-1}(T), x + n-i,k_{i+1}-\delta_{\nwarrow}(\mu_{i+1})+c_{i}(T),\ldots).$$
Let us define a sijection $\Lambda_{n,i} \colon \AR_n \Rightarrow (\{\cdot\},\{\})$. For $\mu' = (\nwarrow,\ldots,\nwarrow,\nearrow,\ldots,\nearrow)$ (with $i$ $\nwarrow$'s), take $\Lambda_{n,i}(\mu') = \cdot$ and $\Lambda_{n,i}(\cdot) = \mu'$. For every other $\mu$, take the smallest $p$ so that $\mu_p \neq \mu'_p$. If $p \leq i$, replace $\nearrow$ with $\nenwarrow$ and vice versa in position $p$ to get $\Lambda_{n,i}(\mu)$ from $\mu$, and if $p > i$, replace $\nwarrow$ with $\nenwarrow$ and vice versa in position $p$ to get $\Lambda_{n,i}(\mu)$ from $\mu$. If $\mu \neq \mu'$, $(\delta_{\nearrow}(\mu_1),\ldots,\delta_{\nearrow}(\mu_{i-1}),\delta_{\nwarrow}(\mu_{i+1}),\ldots,\delta_{\nwarrow}(\mu_n))$ are unaffected by this sijection, so it induces a sijection
\begin{multline*}
  \bigsqcup_{\mu \in \AR_n} \GT(\ldots,k_{i-1}+\delta_{\nearrow}(\mu_{i-1})+c_{i-1}(T), x + n-i,k_{i+1}-\delta_{\nwarrow}(\mu_{i+1})+c_{i}(T),\ldots) \\
  \Rightarrow \GT(\ldots,k_{i-1}+c_{i-1}(T), x + n-i,k_{i+1}+c_{i}(T),\ldots).
\end{multline*}
We switch disjoint unions again, and we get
$$\bigsqcup_{i=1}^n \bigsqcup_{T \in \AP_{n}}  \GT(k_1+c_1(T),\ldots,k_{i-1}+c_{i-1}(T), x + n-i,k_{i+1}+c_{i}(T),\ldots,k_{n} + c_{n-1}(T)).$$
For chosen $i$, use Proposition \ref{prop:sijections} (3) for $\psi = \Psi_{n,i} \circ \Psi_{n,n}^{-1}$ and $\varphi_t = \id$. We get a sijection to
$$\bigsqcup_{i=1}^n \bigsqcup_{T \in \AP_{n}}  \GT(k_1+c_1(T),\ldots,k_{i-1}+c_{i-1}(T), x + n-i,k_{i+1}+c_{i+1}(T),\ldots,k_{n} + c_{n}(T)).$$
If we switch disjoint unions one last time, we can use the sijection $\tau^{-1}$ (see Problem \ref{prob:tau}), and we get
$$\bigsqcup_{T \in \AP_n} \GT(d(\q k,T)) = \SGT(\q k).$$
This completes the construction of $\Phi_4$ and therefore of $\Phi$.
\end{proof}

\begin{problem}
 Given $\q k = (k_1,\ldots,k_n) \in \Z^n$ and $x \in \Z$, construct a sijection
 $$\Gamma = \Gamma_{\q k,x} \colon \MT(\q k) \Rightarrow \SGT(\q k).$$
\end{problem}
\begin{proof}[Construction]
 The proof is by induction on $n$. For $n = 1$, both sides consist of one (positive) element, and the sijection is obvious. Once we have constructed $\Gamma$ for all lists of length less than $n$, we can construct $\Gamma_{\q k,x}$ as the composition of sijections
 $$\MT(\q k) \stackrel{\Xi_{\q k}}{\Longrightarrow} \bigsqcup_{\mu \in \AR_n} \bigsqcup_{\q l \in e(\q k,\mu)} \MT(\q l) \stackrel{\sqcup \sqcup \Gamma}{\Longrightarrow} \bigsqcup_{\mu \in \AR_n} \bigsqcup_{\q l \in e(\q k,\mu)} \SGT(\q l) \stackrel{\Phi_{\q k,x}}{\Longrightarrow} \SGT(\q k),$$
 where $\sqcup \sqcup \Gamma$ means $\bigsqcup_{\mu \in \AR_n} \bigsqcup_{\q l \in e(\q k,\mu)} \Gamma_{\q l,x}.$
\end{proof}

Running the code shows that the main sijection $\Gamma$ indeed depends on the choice $x$. As an example, take $\q k = (1,2,3)$. In this case, $\MT(\q k)$ has $7$ positive elements, and $\SGT(\q k)$ has $10$ positive and $3$ negative elements. For $x = 0$, the sijection is given by
$$\mtthree 1 1 2 1 2 3 \leftrightarrow \left(\mtthree 1 1 1 1 1 1,\mttwo  \searrow \searrow \searrow\right) \qquad \mtthree 2 1 2 1 2 3 \leftrightarrow \left(\mtthree 2 1 2 1 2 2,\mttwo  \searrow \searrow \swarrow\right) \qquad \mtthree 1 1 3 1 2 3 \leftrightarrow \left(\mtthree 1 1 2 1 2 2,\mttwo  \searrow \searrow \swarrow\right)$$
$$\mtthree 2 1 3 1 2 3 \leftrightarrow \left(\mtthree 2 2 3 2 2 3,\mttwo  \swarrow \searrow \swarrow\right) \qquad \mtthree 3 1 3 1 2 3 \leftrightarrow \left(\mtthree 3 2 3 2 2 3, \mttwo  \swarrow \searrow \swarrow\right) \qquad \mtthree 2 2 3 1 2 3 \leftrightarrow \left(\mtthree 2 2 2 3 1 2, \mttwo  \swarrow \seswarrow \searrow\right)$$
$$\mtthree 3 2 3 1 2 3 \leftrightarrow \left(\mtthree 3 3 3 3 3 3,\mttwo  \swarrow \swarrow \swarrow\right) \qquad \left(\mtthree 2 2 2 2 2 3, \mttwo \swarrow \searrow \swarrow\right) \leftrightarrow \left(\mtthree 2 2 2 2 2 2,\mttwo  \seswarrow \searrow \swarrow\right)$$
$$\left(\mtthree 2 2 2 2 3 1,\mttwo \searrow \swarrow \seswarrow\right) \leftrightarrow \left(\mtthree 2 2 2 2 2 2,\mttwo  \swarrow \searrow \seswarrow\right) \qquad \left(\mtthree 2 2 2 1 2 2,\mttwo \searrow \searrow \swarrow\right) \leftrightarrow \left(\mtthree 2 2 2 2 2 2,\mttwo  \searrow \seswarrow \swarrow\right)$$
while for $x = 1$, it is given by
$$\mtthree 1 1 2 1 2 3 \leftrightarrow \left(\mtthree 1 1 1 1 1 1,\mttwo  \searrow \searrow \searrow\right) \qquad \mtthree 2 1 2 1 2 3 \leftrightarrow \left(\mtthree 2 2 2 2 2 3,\mttwo \swarrow \searrow \swarrow\right) \qquad \mtthree 1 1 3 1 2 3 \leftrightarrow \left(\mtthree 1 1 2 1 2 2,\mttwo  \searrow \searrow \swarrow\right)$$
$$\mtthree 2 1 3 1 2 3 \leftrightarrow \left(\mtthree 2 2 3 2 2 3,\mttwo  \swarrow \searrow \swarrow\right) \qquad \mtthree 3 1 3 1 2 3 \leftrightarrow \left(\mtthree 3 2 3 2 2 3,\mttwo  \swarrow \searrow \swarrow\right) \qquad \mtthree 2 2 3 1 2 3 \leftrightarrow \left(\mtthree 2 2 2 1 2 2,\mttwo \searrow \searrow \swarrow\right)$$
$$\mtthree 3 2 3 1 2 3 \leftrightarrow \left(\mtthree 3 3 3 3 3 3,\mttwo  \swarrow \swarrow \swarrow\right) \qquad \left(\mtthree 2 1 2 1 2 2,\mttwo  \searrow \searrow \swarrow\right) \leftrightarrow \left(\mtthree 2 2 2 2 2 2,\mttwo  \seswarrow \searrow \swarrow\right)$$
$$\left(\mtthree 2 2 2 3 1 2, \mttwo  \swarrow \seswarrow \searrow \right) \leftrightarrow \left(\mtthree 2 2 2 2 2 2, \mttwo  \swarrow \searrow \seswarrow \right) \qquad \left(\mtthree 2 2 2 2 3 1,\mttwo \searrow \swarrow \seswarrow \right) \leftrightarrow \left(\mtthree 2 2 2 2 2 2, \mttwo  \searrow \seswarrow \swarrow\right)$$

\section{Concluding remarks}

\subsection*{Future work}

In this article, we have presented the first bijective proof of the operator formula. The operator formula is the main tool for non-combinatorial proofs of several results where alternating sign matrix objects are related to plane partition objects, or simply for showing that $n \times n$ ASMs are enumerated by \eqref{asm}.

\begin{itemize}
\item The operator formula was used in \cite{Fis07} to show that $n \times n$ ASMs are counted by \eqref{asm} and, more generally, to count ASMs with respect to the position of the unique $1$ in the top row.
\item While working on this project, we actually realized that the final calculation in \cite{Fis07} also implies that ASMs are equinumerous with DPPs without having to use Andrews' result \cite{And79} on the number of DPPs; more generally, we can even obtain the equivalence of the refined count of $n \times n$ ASMs with respect to the position of the unique $1$ in the top row and the refined count of DPPs with parts no greater than $n$ with respect to the number of parts equal to $n$. This was conjectured in  \cite{MilRobRum83} and first proved in \cite{BehDifZin12}.
\item In \cite{Fis19a}, the operator formula was used to show that ASTs with $n$ rows are equinumerous with TSSCPPs in a $2n \times 2n \times 2n$-box. Again we do not rely on Andrews' result \cite{And94} on the number of TSSCPPs and we were actually able to deal with a refined count again (which has also the same distribution as the position of the unique $1$ in the top row of an ASM).
\item In \cite{Fis19b}, we have considered alternating sign trapezoids (which generalize ASTs) and, using the operator formula, we have shown that they are equinumerous with objects generalizing DPPs. These objects were already known to Andrews and he actually enumerated them in \cite{And79}. Later Krattenthaler \cite{Kra06} realized that these more general objects are (almost trivially) equivalent to cyclically symmetric lozenge tilings of a hexagon with a triangular hole in the center. Again we do not rely on Andrews' enumeration of these generalized DPPs, and in this case we were able to include three statistics.
\end{itemize}

We plan to work on converting the proofs just mentioned into bijective proofs. For those mentioned in the first and second bullet point, this has already been worked out.
The attentive reader will have noticed that working out all of them will link all four known classes of objects that are enumerated by \eqref{asm}.

\subsection*{Computer code}

As mentioned before, we consider computer code for the constructed sijections an essential part of this project. The code (in python) is available at 
\begin{center}
\url{https://www.fmf.uni-lj.si/~konvalinka/asmcode.html}.
\end{center}
 All the constructed sijections are quite efficient.  If run with pypy, checking that $\Gamma_{(1,2,3,4,5),0}$ is a sijection between $\MT(1,2,3,4,5)$ (with $429$ positive elements and no negative elements) and $\SGT(1,2,3,4,5)$ (with $18913$ positive elements and $18484$ negative elements) takes less than a minute. Of course, the sets involved can be huge, so checking that  $\Gamma_{(1,2,3,4,5,6),0}$ is a sijection between $\MT(1,2,3,4,5,6)$ (with $7436$ positive elements and no negative elements) and $\SGT(1,2,3,4,5,6)$ (with $11167588$ positive elements and $11160152$ negative elements) took almost 20 hours.

\end{document}